\documentclass{article}[11pt]
\usepackage{graphicx}

\usepackage{amsmath,amsfonts,amssymb,latexsym,epsfig}
\usepackage{mathrsfs}
\usepackage{verbatim}
\usepackage{latexsym}
\usepackage{amsthm}
\usepackage{amssymb}
\usepackage{graphics}
\usepackage{amsbsy}
\usepackage{fullpage}
\usepackage{enumerate}
\usepackage{times}
\usepackage{multirow}
\usepackage{soul}
\usepackage[normalem]{ulem}

\usepackage{pgfplots}
\usepackage{tikz}
\usetikzlibrary{decorations.pathreplacing,decorations.markings,arrows}
\newcommand{\R}{\mathbb{R}}

\newcommand{\T}{\textsf{T}}

\newtheorem{theorem}{Theorem}[section]
\newtheorem{lemma}[theorem]{Lemma}
\newtheorem{remark}[theorem]{Remark}
\usepackage{caption}
\usepackage{subcaption}
\numberwithin{equation}{section}
\usepackage{amscd}
\usepackage{tikz}
\usepackage{xcolor}
\usepackage[english]{babel}
\usepackage[latin1]{inputenc}
\usepackage{times}
\usepackage[T1]{fontenc}
\usepackage{graphicx}
\usepackage{amsmath,amssymb}
\usepackage{slidesec}
\usepackage{verbatim}
\usepackage{url}
\usepackage{epsf}
\usepackage{amsmath}
\usepackage{graphics}
\usepackage{amsfonts}
\usepackage{amsbsy}
\usepackage{lscape}
\usepackage{enumerate}
\usepackage{amsthm}
\usepackage{amssymb}
\usepackage{exscale}
\usepackage{algorithm}
\usepackage{algorithmic}

\newcommand\norm[1]{\lVert#1\rVert}
\newcommand{\modk}[1]{\textcolor{black}{#1}\index {#1}}

\newcommand{\modkz}[1]{\textcolor{black}{#1}\index {#1}}
\newcommand{\modp}[1]{\textcolor{black}{#1}\index {#1}}
\newcommand{\modg}[1]{\textcolor{black}{#1}\index {#1}}

\begin{document}
\title{On the connections between optimization algorithms, Lyapunov functions, and differential equations: Theory and insights}

\author{Paul Dobson $^{1,2}$ \and J. M. Sanz-Serna $^{3}$ \and Konstantinos C. Zygalakis  $^{2,4}$}

\maketitle

\begin{abstract}
We revisit the general framework introduced by Fazylab \emph{et al.} (SIAM J. Optim. 28, 2018) to construct Lyapunov functions for optimization algorithms in discrete and continuous time. For smooth, strongly convex objective functions, we relax the requirements necessary for such a construction. As a result we are able to prove for Polyak's ordinary differential equations and for a two-parameter family of Nesterov algorithms rates of convergence that improve on those available in the literature. We analyse the interpretation of Nesterov algorithms as discretizations of the Polyak equation. We show that the algorithms are instances of Additive Runge-Kutta integrators and discuss the reasons why most discretizations of the differential equation do not result in optimization algorithms with acceleration.
We also introduce a modification of Polyak's equation and study its convergence properties.
Finally we extend the general framework to the stochastic scenario and consider an application to random algorithms with acceleration for overparameterized models; again we are able to prove convergence rates that improve on those in the literature.
\end{abstract}

\footnotetext[1]{School of Mathematical and Computer Sciences, Heriot-Watt University, Edinburgh, EH144AS, Scotland,UK.}
\footnotetext[2]{Maxwell Institute for Mathematical Sciences, The Bayes Centre, 47 Potterrow, EH8 9BT, Edinburgh, Scotland, UK.}
\footnotetext[3]{Departamento de Matem\'aticas,
	Universidad Carlos III de Madrid,
 	Legan\'es (Madrid), Spain}
\footnotetext[4]{School of Mathematics, University of Edinburgh, James Clerk Maxwell Building, Edinburgh, EH9 3FD, Scotland, UK}

\maketitle

\section{Introduction}
In this paper we contribute to the literature that explores the relations between optimization algorithms, differential equations and Lyapunov functions \cite{LO20,PS17,SSKZ21,WRJ21}.
As is well known, in order to find a minimizer $x^\star$ of a differentiable function
$f: \R^{d} \rightarrow \R$, the simplest technique is given by the gradient descent (GD) algorithm
\begin{equation} \label{eq:gd}
x_{k+1}=x_{k}-\alpha \nabla f(x_{k}),
\end{equation}
which can be seen as the result of discretizing the gradient flow (GF) ordinary differential equation (ODE)
\begin{equation}\label{eq:gf}
\frac{dx(t)}{dt}=-\nabla f(x(t)) %\quad x(0)=x_{0}  \in \R^{d}
\end{equation}
by means of Euler's rule, the simplest conceivable integrator. While, under very general hypotheses
($f$ bounded from below and $\nabla f$ Lipschitz), the iterates \eqref{eq:gd} will converge to a stationary point of $f$ if $\alpha$ is suitably chosen,  it is standard \cite{N14} to analyze GD when the attention is restricted to functions $f$ that possess additional properties.
For appropriate choices of  $\alpha$,  $f(x_{k})-f(x^\star)$ converge at a rate $\mathcal{O}(1/k)$ when $f \in \mathcal{F}_{L}$ (the set of convex functions with $L$-Lipschitz gradient), while when $f \in \mathcal{F}_{m,L}$ (the set of $m$-strongly convex functions with $L$-Lipschitz gradient) one can show that  $f(x_{k})-f(x^{\star})$  converge with rate  $\mathcal{O}(((\kappa-1)/(\kappa+1))^{k})$, where $\kappa$ denotes the condition number $\kappa=L/m$.

It is of course possible to improve on the rates provided by GD while staying with first-order information, i.e.\ without resorting to information on higher derivatives of $f$. For instance, the celebrated Nesterov's algorithm
\begin{subequations} \label{eq:nest}
\begin{eqnarray}
x_{k+1}&=&x_{k}-\alpha_{k}\nabla f(y_{k}), \\
y_{k+1}&=&x_{k+1}+\beta_{k}(x_{k+1}-x_{k})
\end{eqnarray}
\end{subequations}
converges with rate $\mathcal{O}(1/k^{2})$ for $f \in \mathcal{F}_{L}$ and with rate $\mathcal{O}\big( ((\sqrt{\kappa}-1)/(\sqrt{\kappa}+1) )^{k} \big)$ when $f \in \mathcal{F}_{m,L}$,  for appropriate choices of $\alpha_{k}$, $\beta_k$ (which depend on the class $\mathcal F$ of functions under consideration). This improvement in convergence rate is known as acceleration. The rates quoted for \eqref{eq:nest} are  nearly optimal in terms of what a first-order algorithm can achieve  for both classes of functions \cite{N14}.

As is the case for GD, Nesterov algorithm is related to ODEs, even though the connection was not mentioned in the original paper \cite{N83}. The well-known contribution \cite{SBC16} showed, that,
when $\alpha_k$, $\beta_k$  are tailored for $f \in \mathcal{F}_{L}$, \eqref{eq:nest} provides a numerical discretization of
\begin{equation*}
\ddot{x}(t)+\frac{r}{t}\dot{x}(t)+\nabla f(x(t))=0.
\end{equation*}

\modp{Turning now to the case of $f\in \mathcal{F}_{m,L}$, it is customary to choose, \modg{independently of \(k\)}, $\alpha_k=\alpha$ and $\beta_k=\beta$. For suitable choices of $\alpha$ and $\beta$ the algorithm \eqref{eq:nest} can be seen as a sophisticated (see e.g.\ \cite{SSKZ21}) numerical discretization of the time-independent ODE}
%\modkz{In the case  of $f \in \mathcal{F}_{m,L}$ if one chooses $\alpha_k=\alpha$, $\beta_k=\beta$  for suitable choices of $\alpha,\beta$}
%\eqref{eq:nest} can be seen as a sophisticated (see e.g.\ \cite{SSKZ21}) numerical discretization of the ODE
\begin{equation}\label{eq:polyak}
\ddot{x}(t)+\bar{b} \sqrt{m}\dot{x}(t)+\nabla f(x(t))=0
\end{equation}
considered by Polyak \cite{P64}.\footnote{Following \cite{SSKZ21}, we use overbars for parameters associated to ODEs.} Polyak showed that the heavy ball algorithm, a straightforward discretization of
\eqref{eq:polyak}, exhibits acceleration when applied to \emph{quadratic} $f$.

 The connection between differential equations and optimization algorithms,  further highlighted in   \cite{SRB17}, has led to a, by now large, number of research works that proposed accelerated algorithms both in  Euclidean and non-Euclidean geometry,  based on discretizations of second-order
dissipative ODEs (see e.g.\ \cite{WWJ16,KBB15}). Furthermore, the links with Hamiltonian dynamics have  motivated contributions that  construct or interpret optimization algorithms using concepts such as shadowing \cite{OL19}, symplecticity
\cite{BJW18,BDF19,MJ19,MJ20,SDS19}, discrete gradients \cite{ERRS18},  or backward error analysis \cite{FJV21}. A common element of the analysis presented in many of these papers  is the construction of   discrete Lyapunov functions that were used to investigate the convergence rate of the optimization algorithms.
The reference \cite{FRMP18}, based on the control theoretic view of optimization algorithms suggested in  \cite{LRP16}, has given a general methodology to find convergence rates by means of Lyapunov functions. Applications of this technique may be seen in \cite{SSKZ21}.

In this work, we restrict our attention to the case of strongly convex functions
and modify the general control theory  framework in \cite{FRMP18,SSKZ21}.
\modg{Using the modified framework,} we relax some of the conditions needed to obtain a
Lyapunov function both in  continuous and  discrete time \modkz{(see Theorems \ref{theo:main_con} and  \ref{theo:main_disc}}). {In the new framework, we construct a  Lyapunov function for  \eqref{eq:polyak} \modkz{that allows us to prove in Theorem \ref{thm:polyakconvergence}},  for each choice of the friction parameter $\bar b$, a convergence rate that improves on the rate established in \cite{SSKZ21}. We show that, for $f \in \mathcal{F}_{m,L}$, $\bar b$ may be chosen to guarantee rates arbitrarily close to $\sqrt{2m}$; this is to be compared with the best rate
$\sqrt{m}$ that may be proved in the approach of \cite{SSKZ21,FRMP18}. \modk{Furthermore, this analysis closes the gap between  quadratic and non-quadratic objective functions; in particular, for $\bar{b}>3\sqrt{2}/2$ the convergence rate given in this paper is equal to the rate for quadratic objective functions, showing in this case that the rate is sharp.} Similarly, in the discrete time setting, we obtain a new Lyapunov function for a two-parameter family of Nesterov optimization methods \eqref{eq:nest}. This allows us to
prove \modkz{in Theorem \ref{thm:nestconvergence}}, for a suitable choice of parameters, \modkz{a  convergence rate for $\norm{x_{k}-x^{\star}}^{2}$ arbitrarily close to $ 1 -\sqrt{2}/\sqrt{\kappa}$,  an improvement over the best rate
$1-1/\sqrt{\kappa}$ available in the literature \cite{N14}. Notably, this improved rate is achieved for a  value of the parameter $\beta$ different from the one previously suggested \cite{N14}.}

In addition, the modified  framework is
\begin{enumerate}
\item Used to study a perturbation of the GF equation \eqref{eq:gf} that leads to a new second-order ODE related to \eqref{eq:polyak}. Discretizations of this ODE have the potential of yielding optimization algorithms with acceleration.
\item Extended to account for  stochastic optimization algorithms. This extension is illustrated in the case of accelerated algorithms for over-parameterized models, where again we are able to prove rates  better than those available in the literature \cite{VBS19}.
\end{enumerate}

A final contribution of this work is to interpret \eqref{eq:nest} as a member of the class of additive Runge-Kutta methods \cite{Cooper80} and explain the (rather demanding) structural conditions that discretizations of \eqref{eq:polyak} should satisfy in order to lead to accelerated algorithms.}

The rest of the paper is organized as follows. In Section \ref{sec:pre} we describe the control theoretic framework both in continuous and in discrete time
and formulate general results for the construction of Lyapunov functions. We then in Section \ref{sec:polyak} study the convergence properties of \eqref{eq:polyak} as well as the family of algorithms \eqref{eq:nest}. Section \ref{sec:opt} analyses the connections between  algorithms of the form \eqref{eq:nest} and the ODE \eqref{eq:polyak}.
%We  highlight that the algorithms  may be understood as additive Runge-Kutta discretizations of the ODE, and comment on the  structural conditions that discretizations of \eqref{eq:polyak} need to satisfy  to achieve acceleration.
In Section \ref{sec:polyakplus}  we study a perturbation of the GF ODE. Finally in Section \ref{sec:over} we extend our approach to  stochastic optimization algorithms and in particular consider accelerated algorithms for over-parameterized models.
\section{Preliminaries}
\label{sec:pre}

\subsection{Control theoretic formulation}

We start by discussing a control theoretical formulation
\cite{LRP16,FRMP18} of optimization algorithms both in continuous and in discrete time.

In the continuous time setting, we will consider the following format
\begin{equation} \label{eq:con_system}
\dot{\xi}(t)=\bar{A}\xi(t)+\bar{B}u(t), \quad x(t)=\bar{C}\xi(t), \quad u(t)=\nabla f(x(t)), \quad  t\geq0,
\end{equation}
where  $\xi(t) \in \R^{n}$ is the state, $x(t) \in \R^{d} (d \leq n)$ the feedback output mapped to the input $u(t)=\nabla f(x(t))$. Fixed points of \eqref{eq:con_system} satisfy
\[
0=\bar{A}\xi^\star, \quad x^\star=\bar{C}\xi^\star, \quad u^\star=\nabla f(x^\star);
\]
in the optimization context \(u^\star = 0\) and \(x^\star\) is the minimizer we seek. Both the GF equation \eqref{eq:gf} and Polyak's ODE \eqref{eq:polyak} can be cast in the format \eqref{eq:con_system}. For GF,  $n=d$, $\xi=x$, and $\bar{A}=0_{d}$ , \modkz{$ \bar{B}=-I_{d}$}, $\bar{C}=I_{d}$,
while for \eqref{eq:polyak}  $n=2d$,  \(\xi = [\dot{x}^\T,x^\T]^\T\), and
\[
\bar{A} = \left[\begin{matrix}-\bar{b}\sqrt{m} I_{d} & 0_{d}\\ I_{d} & 0_{d}\end{matrix}\right],\quad
\bar{B} = \left[\begin{matrix}- I_{d} \\  0_{d}\end{matrix}\right],\quad
\bar{C} = \left[\begin{matrix}0_{d}& I_{d}\end{matrix}\right].
\]

In the discrete-time setting, we consider the  formulation
\begin{subequations} \label{eq:control_disc}
\begin{align}
\xi_{k+1} &=A\xi_{k}+Bu_{k}, \\
u_{k} & = \nabla f(y_{k}), \\
y_{k} & =C\xi_{k}, \\
x_{k} &=E \xi_{k},
\end{align}
\end{subequations}
where $\xi_{k} \in \R^{n}$ is the state, $u_{k} \in \R^{d}$ is the input
$(d \leq n)$, $y_{k}  \in \R^{d}$ is the feedback output that is mapped  to
$u_{k}$ by the nonlinear map $\nabla f$, \modkz{and $x_{k}$ is the output used for evaluating the objective function. Note that this differs from the continuous setting (see equation~\eqref{eq:con_system}) where the objective function is evaluated at the feedback output $x(t)$. }
%Here $x_{k}$ is the
%approximation to the mimimizer \(x^\star\).

The algorithms GD  \eqref{eq:gd} and Nesterov's \eqref{eq:nest} in the particular case $\alpha_{k}=\alpha, \beta_{k}=\beta$ we will be focusing on below are easily written in this format. For GD,  $n=d$ and
$
A=0_{d}$, $B=-I_{d}$, $C=I_{d}$, $E=I_{d}$,
while for  \eqref{eq:nest},   $n=2d, \xi_{k}=[x_{k-1}^\T,x_{k}^\T]^{\T}$ and,
\begin{equation}\label{eq:nest_cont_framework}
A=\left[\begin{matrix}
0_{d} & I_{d} \\
-\beta & (\beta+1)I_{d}
\end{matrix} \right], \quad
B=\left[\begin{matrix}
0_{d} \\ -\alpha I_{d}
\end{matrix} \right], \quad
C=\left[\begin{matrix}-\beta I_{d} & (\beta+1)I_{d} \end{matrix} \right], \quad
E=\left[\begin{matrix}
0_{d} & I_{d}.
\end{matrix} \right]
\end{equation}
  The format \eqref{eq:con_system}  can be easily extended  \cite{FRMP18} to cases where $\bar{A}$, $\bar{B}$, \dots depend on $t$. Likewise in \eqref{eq:control_disc} it is possible to let $A$, $B$, \dots depend on $k$. Those extensions are not needed for our purposes here.

\subsection{Matrix inequalities}\label{subsec:MatrixIneq}

 Matrix inequalities may be used to describe different classes of nonlinearities in control theory \cite{MR97}. For the application within optimization see e.g.\ \cite{LRP16,FRMP18}. The key idea here is to express different properties of the function $f$ as matrix inequalities that relate increments in $\nabla f(x)$ and increments in $x$.   For example, a function is $m$-strongly convex if and only if for all $x,y\in\R^d$
\begin{equation*}
    m\lVert x-y\rVert^2\leq (x-y)^\T (\nabla f(x)-\nabla f(y)).
\end{equation*}
This is equivalent to the following matrix inequality: $f$ is $m$-strongly convex if and only if
\begin{equation*} %\label{eq:strong_convex}
    \left[\begin{matrix}
        x-y\\
        \nabla f(x)-\nabla f(y)
    \end{matrix}\right]^\T \left[\begin{matrix}
        -m I_d & \frac{1}{2}I_d\\
        \frac{1}{2}I_d & 0_d
    \end{matrix}\right]\left[\begin{matrix}
        x-y\\
        \nabla f(x)-\nabla f(y)
    \end{matrix}\right] \geq 0.
\end{equation*}
In this work, we will use two additional inequalities for $f\in \mathcal{F}_{m,L}$. If $\nabla f$ is $L$-Lipschitz,  we have
\[
f(x)-f(y) \leq \nabla f(y)^\T (x-y) +\frac{L}{2}\norm{x-y}^{2}
\]
which can be expressed as

\begin{equation} \label{eq:L_smooth}
	f(x)-f(y) \leq  \left[\begin{matrix}
		x-y\\
	\nabla f(y)
	\end{matrix}\right]^\T \left[\begin{matrix}
		\frac{L}{2}I_d & \frac{1}{2}I_d\\
		\frac{1}{2}I_d & 0
	\end{matrix}\right]\left[\begin{matrix}
		x-y\\
		\nabla f(y)
	\end{matrix}\right].
\end{equation}
For  $f\in \mathcal{F}_{m,L}$, we have that
\[
\frac{mL}{m+L}\lVert x-y\rVert^2+\frac{1}{m+L}\lVert \nabla f(x)-\nabla f(y)\rVert^2 \leq (\nabla f(x)-\nabla f(y))^\T (x-y)
\]
which gives rise to:
\begin{equation} \label{eq:L_m_together}
    \left[\begin{matrix}
        x-y\\
        \nabla f(x)-\nabla f(y)
    \end{matrix}\right]^\T \left[\begin{matrix}
        -\frac{mL}{m+L}I_d & \frac{1}{2}I_d\\
        \frac{1}{2}I_d & \frac{-1}{m+L}I_d
    \end{matrix}\right]\left[\begin{matrix}
        x-y\\
        \nabla f(x)-\nabla f(y)
    \end{matrix}\right] \geq 0.
\end{equation}

\subsection{Lyapunov functions for ODEs and their discretizations} \label{subsec:lyap}
A way to study the convergence of the continuous dynamics \eqref{eq:con_system} and their discrete counterparts \eqref{eq:control_disc}  is by using   Lyapunov functions. In the case of continuous dynamics, the references \cite{FRMP18, SSKZ21} use  Lyapunov functions of the form
\begin{equation} \label{eq:cont_lyap}
V(\xi(t),t)=e^{\lambda t} \Big(f(x(t))-f(x^{\star})+(\xi(t)-\xi^\star)^{\T}
\bar{P}(\xi(t)-\xi^\star) \Big),
\end{equation}
 where  $\lambda>0$ and  $\bar{P}$ is an $n\times n$ symmetric matrix. If one can show that, for suitably chosen $\lambda$ and $\bar P$,
 $(d/dt){V} \leq 0$ along solutions of \eqref{eq:con_system}, then
 $$
 e^{\lambda t} \Big(f(x(t))-f(x^{\star})+ (\xi(t)-\xi^\star)^{\T}
\bar{P}(\xi(t)-\xi^\star)\Big)
 \leq
  V(\xi(0),0),
 $$
 which, under the additional assumption that $\bar{P}$ is \emph{positive semidefinite}, $\bar{P}\succeq 0$, leads  to the decay estimate
 \[
f(x(t))-f(x^{\star}) \leq e^{-\lambda t}  V(\xi(0),0).
\]

 In this paper, we relax the hypothesis
 $\bar{P}\succeq 0$ in order  to improve the decay rate $\lambda$. We leverage the fact that the attention is restricted to $f \in \mathcal{F}_{m,L}$ and therefore
 \begin{equation} \label{eq:upper_bound}
     \frac{m}{2}\norm{x(t)-x^{\star}}^{2}\leq f(x(t))-f(x^{\star}),
\end{equation}
so that from \eqref{eq:cont_lyap}, using the relation between $\xi$ and $x$ in \eqref{eq:con_system},
$$
  e^{\lambda t} \Big((\xi(t)-\xi^\star)^{\T}
\widetilde{P} (\xi(t)-\xi^\star)\Big)\leq V(\xi(t),t),
$$
where $\modkz{\widetilde{P}=\bar{P}+(m/2)\bar{C}^{\T}\bar{C}}$. Thus, if $V$ decreases along the dynamics,
$$
(\xi(t)-\xi^\star)^{\T}
\widetilde{P} (\xi(t)-\xi^\star) \leq e^{-\lambda t} V(\xi(0),0),
$$
which, after using \eqref{eq:con_system} once more, leads to the following decay estimate for $x$
($\sigma$ denotes the spectrum of eigenvalues):
\modkz{\begin{eqnarray} \label{eq:conv_cont}
\norm{x(t)-x^{\star}}^{2}
&\leq& \max \sigma(\bar{C}^{\T}\bar{C}) \norm{\xi(t)-\xi^{\star}}^2 \nonumber \\
&\leq & \frac{ \max \sigma({\bar C}^{\T}\bar C)}{\min \sigma(\widetilde{P})}\norm{\xi(t)-\xi^{*}}^{2}_{\widetilde{P}} \\
&\leq & \frac{ \max \sigma({\bar C}^{\T}\bar C)}{\min \sigma(\widetilde{P})}e^{-\lambda t} V(\xi(0),0), \nonumber
\end{eqnarray}
}
provided that $\min \sigma(\widetilde{P})>0$, i.e.\ that $\widetilde{P}\succ 0$. \modkz{Here $$\norm{\xi(t)-\xi^{*}}^{2}_{\widetilde{P}}=(\xi(t)-\xi^\star)^{\T}
\widetilde{P} (\xi(t)-\xi^\star).$$}

The following theorem provides conditions that guarantee that the Lyapunov function \eqref{eq:cont_lyap} is indeed decreasing along the trajectories of \eqref{eq:con_system} so that \eqref{eq:conv_cont} holds. The proof, \modp{see Section \ref{subsec:proof_main_con},} is similar to the proof of Theorem  6.4 in \cite{FRMP18} and relies on computing $(d/dt)V$ along the dynamics and using  the relations \eqref{eq:L_smooth} and \eqref{eq:L_m_together}.
\begin{theorem} \label{theo:main_con} Suppose that, for \eqref{eq:con_system}, there exist $\lambda >0$,  $\sigma \geq 0$ and a symmetric matrix $\bar{P}$ with  $\widetilde{P}:=\bar{P}+(m/2)\bar{C}^{\T}\bar{C} \succ 0$,   that satisfy
$$%\begin{equation} \label{eq:LMI_con}
\bar{T}= \bar M^{(0)}+\bar M^{(1)}+\lambda \bar M^{(2)}+\sigma \bar M^{(3)} \preceq 0
$$%\end{equation}
where
\begin{align*}
\bar M^{(0)} &=
\left[ \begin{matrix}
\bar{P}\bar{A}+\bar{A}^{\T}\bar{P}+\lambda \bar{P} & \bar{P}\bar{B} \\
\bar{B} ^{\T}\bar{P}  & 0
\end{matrix}\right], \\
\bar M^{(1)} &= \frac{1}{2}\left[ \begin{matrix}
0 & (\bar{C} \bar{A} )^{\T} \\
\bar{C} \bar{A} & \bar{C} \bar{B} +\bar{B} ^{\T}\bar{C} ^{\T}
\end{matrix}\right], \\
\bar M^{(2)}&=\left[ \begin{matrix}
\bar{C}^{\T} & 0 \\
0 & I_{d}
\end{matrix}\right]\left[ \begin{matrix}
-\frac{m}{2}I_{d} & \frac{1}{2}I_{d} \\
\frac{1}{2}I_{d} & 0
\end{matrix}\right]\left[ \begin{matrix}
\bar{C} & 0 \\
0 & I_{d}
\end{matrix}\right], \\
\bar M^{(3)} &=\left[ \begin{matrix}
\bar{C}^{\T} & 0 \\
0 & I_{d}
\end{matrix}\right]\left[ \begin{matrix}
-\frac{m L}{m+L}I_{d} & \frac{1}{2}I_{d}, \\
\frac{1}{2}I_{d} & -\frac{1}{m+L}I_{d}
\end{matrix}\right]\left[ \begin{matrix}
\bar{C} & 0 \\
0 & I_{d}
\end{matrix}\right].
\end{align*}
Then for $f \in \mathcal{F}_{m,L}$,  $t \geq 0$, and $V$ given by \eqref{eq:cont_lyap}, the decay estimate
\eqref{eq:conv_cont} holds.
\end{theorem}

\begin{remark}The Lipschitz constant $L$ only appears in $\bar T$ through the matrix $\bar{M}^{(3)}$. Therefore if $\sigma = 0$ the theorem holds for arbitrary $m$-strongly convex $f$.\end{remark}

The case of the discrete dynamics \eqref{eq:control_disc} is completely parallel. The Lyapunov functions considered are of the form
\begin{equation} \label{eq:liap_disc}
V_{k}(\xi)=\rho^{-2 k} \left(a_{0}(f(x)-f(x^\star))+(\xi-\xi^\star)^{\T}
P(\xi-\xi^\star) \right),\qquad \rho \in (0,1),
\end{equation}
with $P$ symmetric and $a_{0}>0$. If one can prove that along the discrete dynamics $V_{k+1}(\xi_{k+1}) \leq V_{k}(\xi_{k})$ then, for   $P\succeq 0$, it is easy to show that
\[
 f(x_{k})- f(x^\star) \leq \rho^{2k} \frac{V_{0}(\xi_{0})}{a_{0}}.
\]
In this paper, for $f \in
\mathcal{F}_{m,L}$,  we relax the assumption $P\succeq 0$ by exploiting the bound \eqref{eq:upper_bound}.
The following theorem summarises the
conditions that guarantee that the Lyapunov function decays
along the dynamics \eqref{eq:control_disc} and provides a rate
of convergence of $x_{k}$ towards $x^{\star}$.

\begin{theorem} \label{theo:main_disc}
Suppose that, for \eqref{eq:control_disc}, there exist $a_{0} >0$,  $\rho \in (0,1)$, $\ell >0$,
and a symmetric matrix $P$, with
$\widetilde{P}:=P+(a_{0}m/2)E^{\T}E \succ 0$,
such that
\begin{equation} \label{eq:LMI_disc}
T= M^{(0)}
+a_{0}\rho^{2}M^{(1)}+a_{0}(1-\rho^{2})M^{(2)}+\ell M^{(3)} \preceq 0,
\end{equation}
where
\[
M^{(0)}=\left[ \begin{matrix}
A^{\T}PA-\rho^{2}P & A^\T PB \\
B^{\T}PA & B^{\T} P B
\end{matrix}\right],
\]
and
\[
M^{(1)}=N^{(1)}+N^{(2)}, \quad M^{(2)}=N^{(1)}+N^{(3)},\quad  M^{(3)}=N^{(4)},
\]
with
\begin{align*}
N^{(1)} &= \left[ \begin{matrix}
EA-C & EB \\
0 & I_{d}
\end{matrix}\right]^{\T}\left[ \begin{matrix}
\frac{L}{2}I_{d} & \frac{1}{2}I_{d} \\
\frac{1}{2}I_{d} & 0
\end{matrix}\right]\left[ \begin{matrix}
EA-C & EB \\
0 & I_{d}
\end{matrix}\right],
\\
N^{(2)} &= \left[ \begin{matrix}
C-E & 0 \\
0 & I_{d}
\end{matrix}\right]^{\T}\left[ \begin{matrix}
-\frac{m}{2}I_{d} & \frac{1}{2}I_{d} \\
\frac{1}{2}I_{d} & 0
\end{matrix}\right]\left[ \begin{matrix}
C-E & 0 \\
0 & I_{d}
\end{matrix}\right],
\\
N^{(3)} &=\left[ \begin{matrix}
C^{\T} & 0 \\
0 & I_{d}
\end{matrix}\right]\left[ \begin{matrix}
-\frac{m}{2}I_{d} & \frac{1}{2}I_{d} \\
\frac{1}{2}I_{d} & 0
\end{matrix}\right]\left[ \begin{matrix}
C & 0 \\
0 & I_{d}
\end{matrix}\right],
\\
N^{(4)} &=  \left[ \begin{matrix}
C^{\T} & 0 \\
0 & I_{d}
\end{matrix}\right]\left[ \begin{matrix}
-\frac{m L}{m+L}I_{d} & \frac{1}{2}I_{d} \\
\frac{1}{2}I_{d} &  -\frac{1}{m+L}I_{d}
\end{matrix}\right]\left[ \begin{matrix}
C & 0 \\
0 & I_{d}
\end{matrix}\right].
\end{align*}
Then, for $f \in \mathcal{F}_{m,L}$, with $V$ given by \eqref{eq:liap_disc}, the sequence $\{x_{k}\}$ satisfies
\begin{equation} \label{eq:conv_disc}
\norm{x_{k}-x^{\star}}^{2} \leq \max \sigma(E^{\T}E) \norm{\xi_k-\xi^{\star}}^2 \leq \frac{ \max \sigma(E^{\T}E)}{\min \sigma(\widetilde{P})} V_{0}(\xi_0)\rho^{2 k}.
\end{equation}
\end{theorem}
\modp{The proof of this theorem can be found in Section \ref{subsec:proof_main_disc}.}

\section{Analysis of Polyak equation and Nesterov's algorithm}
\label{sec:polyak}
We will now use the framework in Section \ref{subsec:lyap} to study the convergence properties of \eqref{eq:polyak}. We will then present an analysis for the convergence properties of the family of algorithms \eqref{eq:nest}. Both analyses will be connected in Section \ref{sec:opt} by means of the theory of numerical methods for ODEs.

\subsection{Continuous time analysis}
\label{subsec:polyak_ODE}
By introducing the variable $v = (1/\sqrt{m}) \dot x$, equation \eqref{eq:polyak} can be rewritten as the system
\modkz{
\begin{subequations}\label{eq:ode1}
\begin{align}
\dot v(t) &= -\bar b \sqrt{m} v(t) -\frac{1}{\sqrt{m}} \nabla f(x(t)),\\
\dot x(t) &= \sqrt{m} v(t).
\end{align}
\end{subequations}
}
The friction parameter $\bar b$ is nondimensional, i.e.\ it does not change when in \eqref{eq:polyak} $t$, $x$ or $f$ are rescaled.
The scaling factor $\sqrt{m}$ in the definition of $v$ has been introduced to ensure that  \modkz{$v(t)$ shares the dimensions\footnote{\modp{Note that one could rescale time $t\mapsto t/\sqrt{m}$ to obtain the system \eqref{eq:ode1} with $m=1$ and $f$ replaced by $f/m$, which is a $1$-strongly convex function. It would then be   sufficient to establish convergence for the case $m=1$; and all other choices of $m$ would be recovered via this scaling.}} of $x(t)$. If we now  set \(\xi(t) = [v^{\T}(t),x^{\T}(t)]^{\T}\)}, then \eqref{eq:ode1} is of
the form \eqref{eq:con_system} with
\begin{equation} \label{eq:polyak1}
\bar A = \left[\begin{matrix}-\bar b\sqrt{m}I_d & 0_d\\ \sqrt{m} I_d& 0_d\end{matrix}\right],\quad
 \bar B = \left[\begin{matrix}-(1/\sqrt{m}) I_d\\ 0_d\end{matrix}\right],
\quad
 \bar C = \left[\begin{matrix} 0_d & I_d\end{matrix}\right].
\end{equation}

\modp{When using \eqref{eq:ode1} to minimise a given objective function, an estimate of the strong convexity constant $m$ is requiered. Of course, underestimation of the true strong convexity constant will result in a slower convergence of the solution to the fixed point. The following Theorem states the convergence rate towards the minimizer for the solutions of \eqref{eq:polyak1}}

\modp{
	\begin{theorem}\label{thm:polyakconvergence}
		Fix $m>0$, $\overline{b}>0$ with $\overline{b}\neq 3\sqrt{2}/2$ and $f\in \mathcal{F}_{m,L}$ with $L>m$. If $\xi(t)=(v^{T}(t),x^{T}(t))^{T}$ is the solution of \eqref{eq:ode1} with parameter $\overline{b}$ and initial condition $\xi_{0}=(v^{T}_0,x^{T}_0)^{T}$, then
		\begin{equation} \label{eq:cont_polak_convergence}
			\norm{x(t)-x^{\star}}^{2}  \leq C e^{-\bar{r}\sqrt{m}t} \left(f(x_0)-f(x^\star)+\norm{\xi_{0}-\xi^{*}}^{2}_{\tilde{P}}\right),
		\end{equation}
		where
		\begin{eqnarray*}
		\tilde{P}&=&\frac{m}{2}\left[\begin{matrix}
		1 & \bar{r} \\
		\bar{r} & \frac{\bar{r}^{2}}{2}+1
		\end{matrix}\right] \otimes I_{d},\\ C&=&\frac{1}{\min{\sigma(\tilde{P})}}=\frac{8}{m(\bar{r}^{2}+4-\bar{r}\sqrt{\bar{r}^{2}+16})},
		\end{eqnarray*}
		and
		\begin{equation} \label{eq:r_cases}
			\bar{r}=\begin{cases}
				\frac{2\bar{b}}{3}, & \bar{b}<\frac{3\sqrt{2}}{2},\\
				\overline{b}-\sqrt{\bar{b}^2-4}, & \bar{b}>\frac{3\sqrt{2}}{2}.
			\end{cases}
		\end{equation}
		\end{theorem}		
}

\begin{proof}[Proof of Theorem \ref{thm:polyakconvergence}]
\modkz{The proof uses Theorem \ref{theo:main_con} with the following setting (see Remarks \ref{rem:choice1}, \ref{rem:choice2} for additional explanation)
\begin{equation} \label{eq:choices}
\lambda=\sqrt{m}\bar{r}, \quad \bar{P}=\frac{m}{2}\left [ \begin{matrix}1 & \bar{r} \\
		\bar{r} & \frac{\bar{r}^{2}}{2} \end{matrix} \right] \otimes I_{d}, \qquad \sigma=0.
\end{equation}
The matrices $\bar{A}, \bar{B}, \bar{C}$ are as in \eqref{eq:polyak1} and $\bar{T}$ is then given by
\[
\bar{T}=\widehat{\bar{T}} \otimes I_{d}, \qquad
\widehat{\bar{T}}=\frac{m^{3/2}}{2}\left[ \begin{matrix}  3\bar{r}-2\bar{b} & \frac{\bar{r}}{2}\left(3\bar{r}-2\bar{b}  \right) & 0 \\ \frac{\bar{r}}{2}\left(3\bar{r}-2\bar{b}  \right)  & \bar{r}\left(\frac{\bar{r}^{2}}{2}-1 \right)& 0 \\ 0& 0& 0 \end{matrix}\right].
\]
Given the structure of $\bar{T}$, to ensure  $\bar{T} \preceq 0$ we need that the  diagonal elements of $\widehat{\bar{T}}$  be non-positive and the determinant of the leading $2\times 2$ submatrix of $\widehat{\bar{T}}$  be non-negative. Hence, we obtain the conditions
\begin{subequations} \label{eq:T_cond}
\begin{eqnarray}
3\bar{r}-2\bar{b} &\leq& 0, \\
 \bar{r}\left(\frac{\bar{r}^{2}}{2}-1 \right) &\leq& 0, \\
\frac{1}{4}(3\bar{r}-2\bar{b})\bar{r}\left(-\bar{r}^{2}+2\bar{b}\bar{r}-4 \right) &\geq& 0.
\end{eqnarray}
\end{subequations}
Now note that,  since $\bar{r} > 0$,  (\ref{eq:T_cond}b) holds if and only if
$\bar{r} \leq \sqrt{2}$. For, any value of $\bar{r} < \sqrt{2}$, the matrix $\tilde{P}$ is  $\succ 0$.  In addition (\ref{eq:T_cond}a) implies that (\ref{eq:T_cond}c) can only hold if $\bar{r}={2\bar{b}}/{3}$ or if
\begin{equation}\label{eq:quadratic}
-\bar{r}^{2}+2\bar{b}\bar{r}-4=-(\bar{r}-(\bar{b}-\sqrt{\bar{b}^{2}-4}))(\bar{r}-(\bar{b}+\sqrt{\bar{b}^{2}-4})) \leq 0.
\end{equation}
The case $\bar{r}={2\bar{b}}/{3}$ gives a feasible solution for the constraints \eqref{eq:T_cond}, while maintaining the positive definiteness of $\tilde{P}$  for $\bar{b} < {3\sqrt{2}}/{2}$. Observe that in fact this is the maximal solution for $\bar{r}$   in light of (\ref{eq:T_cond}a) giving us thus the first branch of \eqref{eq:r_cases}. We now turn our attention to $\bar{b} \geq {3\sqrt{2}}/{2}$ in which case $\bar{r}-(\bar{b}+\sqrt{\bar{b}^{2}-4}) <0$ for $\bar{r} \leq \sqrt{2}$. Hence by using \eqref{eq:quadratic} it is  now necessary that $\bar{r} \leq \bar{b}-\sqrt{\bar{b}^{2}-4}$. Observe that $\bar{b}-\sqrt{\bar{b}^{2}-4}$ is decreasing in $\bar{b}$ and hence it is sufficient to check if equations (\ref{eq:T_cond}a--b) are satisfied  for $\bar{b}={3\sqrt{2}}/{2}$. This is indeed the case and,  given that we are  interested in the maximal solution for $\bar{r}$, we obtain the second branch of \eqref{eq:r_cases}.
}
\end{proof}
\modg{
\begin{remark} \label{rem:choice1}
Writing the convergence rate $\lambda$ as in \eqref{eq:choices} is just introducing a new parameter $\bar{r}$ that replaces $\lambda$; the new parameter  $\bar{r}$ is non-dimensional and therefore, when rates are reported in terms of $\bar{r}$  as in the figures below, they are independent of the value of $m$.
 On the other hand, similar to what was done in \cite{SSKZ21},  one could have chosen $\sigma>0$ rather than as in \eqref{eq:choices} to obtain a slightly better convergence rate. Nevertheless this improvement in the convergence rate would be  $\mathcal{O}(L^{-1})$ and it would thus be negligible for the ill-conditioned problems of interest.
\end{remark}}

\modkz{\begin{remark}  \label{rem:choice2}
It is perhaps useful to explain how we arrived at the matrix $\bar{P}$ in \eqref{eq:choices}. Firstly we note that the matrix $\bar{A}$ in \eqref{eq:polyak1} is a Kronecker product
\[
\bar{A}= \left[\begin{matrix}-\bar b\sqrt{m}& 0\\ \sqrt{m} & 0\end{matrix}\right] \otimes I_{d}.
\]
The factor $I_{d}$ originates from the dimensionality of $x$ and the $2\times 2$ size of the other factor arises from the fact that \eqref{eq:polyak} is a second order ODE. The matrices $\bar{B},\bar{C}$  have a similar Kronecker product structure. It is thus natural to consider symmetric matrices of the form
\begin{equation} \label{eq:P}
\bar{P}=\widehat{\bar{P}} \otimes I_{d}, \quad \widehat{\bar{P}}=\left[\begin{matrix} \bar{p}_{11} & \bar{p}_{12} \\
\bar{p}_{12} & \bar{p}_{22}
\end{matrix} \right]
\end{equation}
and then $\bar{T}$ will also have a Kronecker product structure
\begin{equation} \label{eq:T}
\bar{T}=\widehat{\bar{T}} \otimes I_{d}, \quad \widehat{\bar{T}}=\left[\begin{matrix}
\bar{t}_{11} & \bar{t}_{12} & \bar{t}_{13} \\
\bar{t}_{12} & \bar{t}_{22} &  \bar{t}_{23} \\
\bar{t}_{13} & \bar{t}_{23} &  \bar{t}_{33}
\end{matrix} \right].
\end{equation}
From \eqref{eq:polyak1}, we find
\begin{align*}
  \bar t_{11} &=  -2\bar b \sqrt{m} \bar p_{11}+2\sqrt{m} \bar p_{12}+\lambda \bar p_{11},\\
  \bar t_{12} &= -\bar b \sqrt{m} \bar p_{12}+\sqrt{m} \bar p_{22}+\lambda \bar p_{12}, \\
  \bar t_{13} &= -(1/\sqrt{m})\bar p_{11}+\sqrt{m}/2, \\
  \bar t_{22} &= \lambda  \bar p_{22}-(m/2)\lambda,\\
  \bar t_{23} &= -(1/\sqrt{m})\bar p_{12}+\lambda/2,\\
  \bar t_{33}&=  0.
\end{align*}
As discussed in Remark \ref{rem:choice1}, we  now set \(\lambda = \sqrt{m}\:{\bar r}\). In addition, since $\bar{t}_{33}=0$, $\widehat{\bar T} \preceq 0$ imposes that  $\bar{t}_{13}=\bar{t}_{23}=0$, which in turn leads to the values of $\bar{p}_{11}, \bar{p}_{12}$ in \eqref{eq:choices}. Now since all the elements  in the third row/column of $\widehat{\bar T}$ vanish, we only have to deal with the leading $2\times 2$ submatrix of $\widehat{\bar T}$.  In particular, it is necessary that   $\Delta(\bar{p}_{22},\bar{r})=\bar t_{11} \bar t_{22} -\bar t^{2}_{12}  \geq 0$. As we  we wish to maximize $\bar r$, we set $(\partial/\partial \bar{p}_{22})\Delta = 0$, giving us the value of $\bar{p}_{22}$ in \eqref{eq:choices}.
Note that the choices of $\bar{p}_{11}, \bar{p}_{12}$ are actually necessary to ensure $\widehat{\bar T} \preceq 0$. On the other hand, while the value of $\bar{p}_{22}$ in \eqref{eq:choices} is not the only one that leads to  $\widehat{\bar T} \preceq 0$, however numerical evidence indicates that is optimal.
\end{remark}
}

\begin{figure}
     \centering
     \begin{subfigure}[b]{0.45\textwidth}
         \centering
         \includegraphics[width=\textwidth]{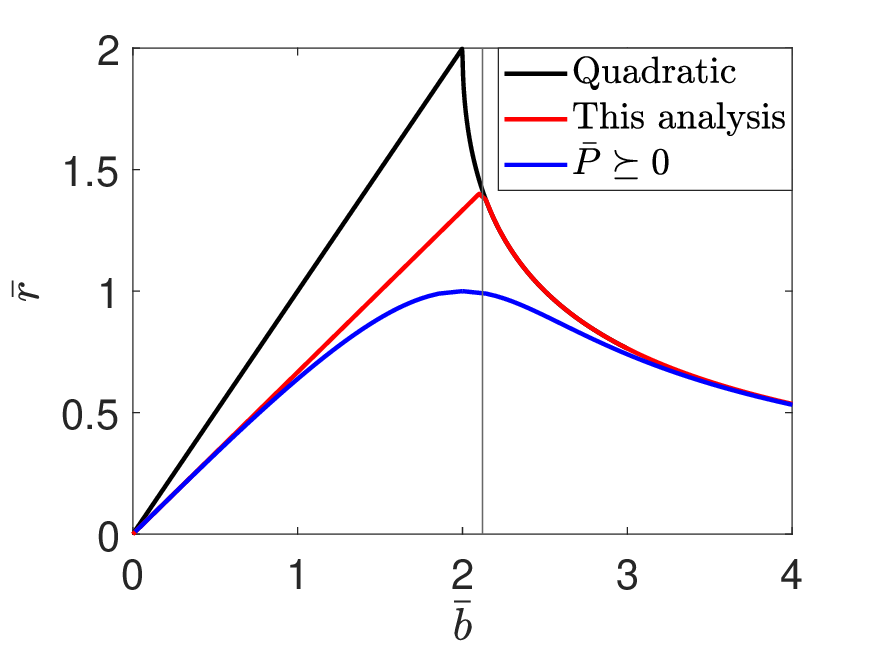}
         \caption{Continuous setting}
         \label{fig:fig1}
     \end{subfigure}
     \hfill
     \begin{subfigure}[b]{0.45\textwidth}
         \centering
         \includegraphics[width=\textwidth]{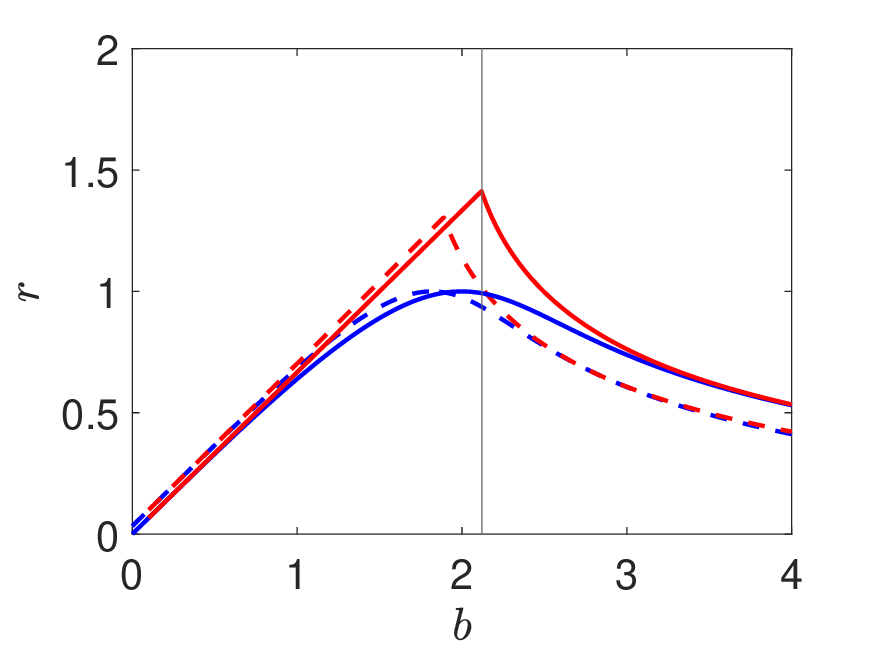}
         \caption{Discrete setting}
         \label{fig:fig2}
     \end{subfigure}
  \caption{The left panel shows the relationship between the rate $\bar{r}=\lambda/\sqrt{m}$ and the parameter $\bar{b}$ in the time-continuous case.
        The right panel shows the relationship between the rate $r$ and the method parameter $b$ in the discrete case when $\delta = \delta_{max}= 1/\sqrt{\kappa}$; the solid curves are for $\kappa = 10^6$ and  the dashed curves are for $\kappa = 10^2$. The red curves correspond to the present analysis and the blue curves correspond to the hypothesis $\bar P\succeq 0$.  The red and blue solid lines on the right are indistinguishable from the red and blue lines on the left.
   }
        \label{fig:conta}
\end{figure}

\modkz{
Theorem \ref{thm:polyakconvergence} implies that it is possible to get all nondimensional rates of convergence $\bar{r}$ in the interval $(0,\sqrt{2})$. Each value of $\bar r\in(0,\sqrt{2})$ may be achieved in two ways, the first by choosing $\bar{b} = 3\bar{r}/2\in (0, 3\sqrt{2}/2)$ and the second by choosing $\bar{b} >3\sqrt{2}/2$. It is not possible to  prove the rate $\bar{r}=\sqrt{2}$ since in this case the matrix $\tilde{P}$ becomes singular. The value of $\bar r$ as a function of  $\bar b$ is represented in Figure~\ref{fig:fig1}, where for comparison we have also provided the best value of $\bar r$ that may be obtained when using the framework in \cite{SSKZ21} that requires $\bar{P}\succeq 0$.  As we can see, the modification of the hypothesis on $\bar P$ allows to prove a significantly better convergence rate.
}

\begin{remark}
If the objective function $f$ is  quadratic, it is of course possible to obtain a sharp bound for the convergence rate $\lambda=\sqrt{m}\bar{r}$ by solving \eqref{eq:polyak} in terms of eigenvalues/vectors.
(See \cite[Section 2.2]{LRP16} for the analysis in the discrete scenario.) Also included  in
 Figure \ref{fig:fig1} is the rate for $m$-strongly convex quadratic problems, which  is maximized for $\bar{b}=2$, where $\lambda = 2\sqrt{m}$. For non-quadratic targets, the rate that may be proved under the hypothesis $\bar P\succeq 0$ in \cite{FRMP18,SSKZ21} is also maximized when $\bar{b}=2$, where $\lambda = \sqrt{m}$. The present analysis proves,  for non-quadratic targets, bounds with rates arbitrarily close to $\lambda = \sqrt{2}\sqrt{m}$, by choosing $\bar b$ close to $3\sqrt{2}/2$. Note that for $\bar b > 3\sqrt{2}/2$ the rate proved here cannot be improved, as it coincides with the rate that the ODE achieves for quadratic objective functions.
 \end{remark}

\subsubsection{A numerical illustration}
We compare numerically the bound provided by the analysis just presented with the corresponding bound when operating within the framework in \cite{SSKZ21}. We use the two-dimensional objective function in
$\mathcal{F}_{m,L}$ given by
\begin{equation} \label{eq:1d}
f(x_1,x_2)=\frac{m}{2}(x_1^{2}+x_2^2)+4(L-m)\log{(1+e^{-x_1})}
\end{equation}
(subindices denote scalar components of the vector $x$) and, for $0\leq t\leq 20$, compute solutions of \eqref{eq:ode1} with a high-order Runge-Kutta algorithm. We report here results for $m =1$, $L=10^6$, when the initial condition is chosen as $x_1(0) = 0$, $x_2(0) = 50$, $v_1(0) = 0$, $v_2(0)=0$.

\begin{figure}[t]
	\centering
	\begin{subfigure}[b]{0.32\textwidth}
		\centering
		\includegraphics[width=\textwidth]{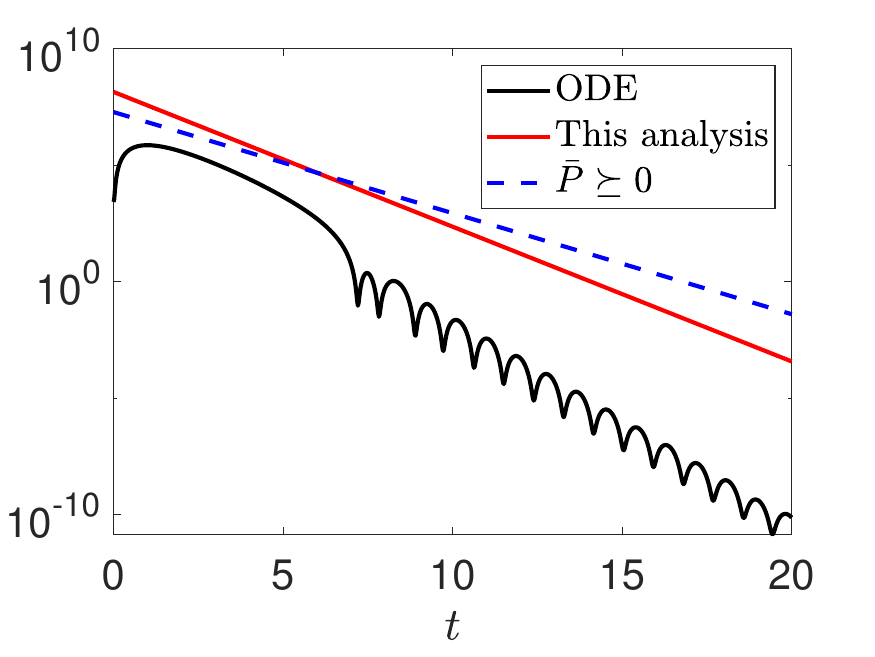}
		\caption{$\bar{b}=2$}
		\label{fig:A}
	\end{subfigure}
	\hfill
	\begin{subfigure}[b]{0.32\textwidth}
		\centering
		\includegraphics[width=\textwidth]{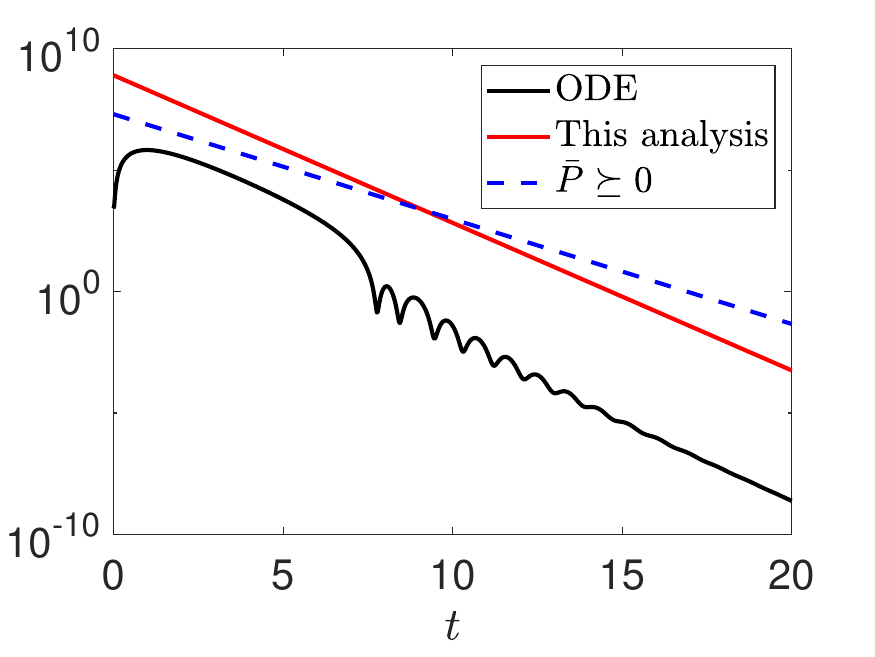}
		\caption{ $\bar{b}=2.1$ }
		\label{fig:B}
	\end{subfigure}
	\hfill
	\begin{subfigure}[b]{0.32\textwidth}
		\centering
		\includegraphics[width=\textwidth]{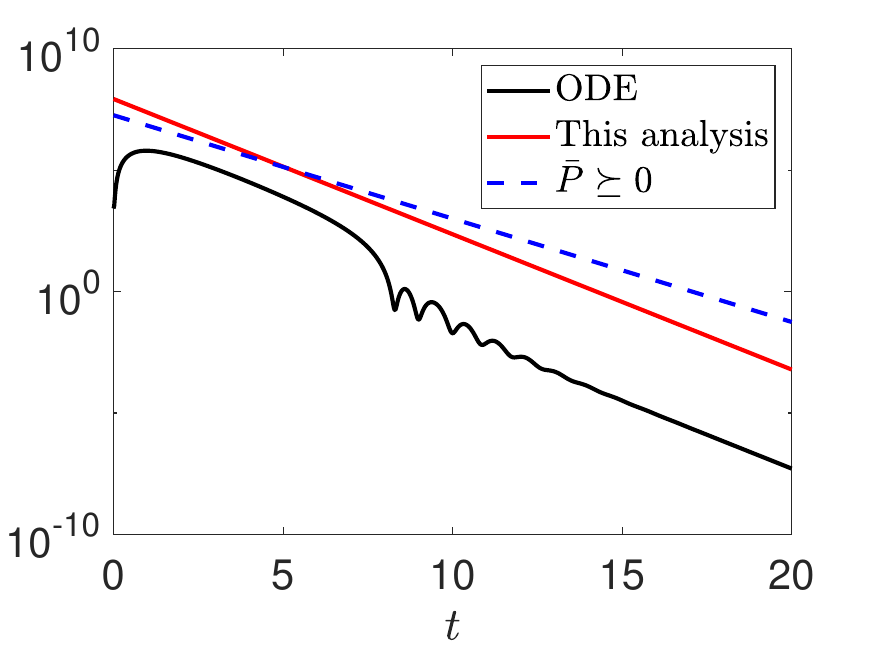}
		\caption{ $\bar{b}=2.2$ }
		\label{fig:C}
	\end{subfigure}
	\caption{Polyak ODE. Bounds for $\norm{x(t)-x^{\star}}^{2}$ for different values of the parameter $\bar{b}$, when $f$ is given by \eqref{eq:1d}, $m=1$, $L=10^{6}$ . }
	\label{fig:cont}
\end{figure}

The first panel in Fig.~\ref{fig:cont} corresponds to $\bar b = 2$, the value that provides the bound with best rate when operating as in \cite{SSKZ21}. The solid straight line gives the bound  \eqref{eq:conv_cont} when $\bar P$ and $\lambda$ are taken as in the analysis in the preceding subsection; one finds $\min\sigma(\tilde P)\approx 0.0195$, $\max\sigma(C^\T C) = 1$ and $\lambda=4/3$. The dashed line gives the bound
\eqref{eq:conv_cont} when $\bar P\succeq 0$ and $\lambda$ are determined as in \cite{SSKZ21}; then
$\min\sigma(\tilde P)\approx 0.1432$, $\max\sigma(C^\T C) = 1$ and $\lambda = 1$. We see how, by relaxing the requirements on $\bar P$, it is possible to prove a larger rate of convergence, at the expense of increasing the factor $1/\min\sigma(\tilde P)$ in \eqref{eq:conv_cont}. In this experiment, the slopes of both straight lines clearly underestimate the true rate of decay in the ODE.

In the central panel of Fig.~\ref{fig:cont}, $\bar b = 2.1$, a value slightly below $3\sqrt{2}/2 \approx 2.1213$. Our analysis yields $\min\sigma(\tilde P)\approx 0.0034$, $\max\sigma(C^\T C) = 1$ and $\lambda=1.400$, while when working as in \cite{SSKZ21}, we get $\min\sigma(\tilde P)\approx 0.1355$, $\max\sigma(C^\T C) = 1$ and the quite pessimistic value $\lambda\approx 0.9950$.

In the final panel, $b=2.2>3\sqrt{2}/2$. Now our analysis has $\min\sigma(\tilde P)\approx 0.0319$, $\max\sigma(C^\T C) = 1$ and $\lambda=1.2835$, and, under the hypotheses of \cite{SSKZ21},
$\min\sigma(\tilde P)\approx 0.1493$, $\max\sigma(C^\T C) = 1$ and  $\lambda\approx 0.9807$. The slope of the the continuous line describes very well the decay behaviour of the ODE (note that, for this value of $\bar b$, the rate proved here cannot be improved, as it coincides with the rate the ODE achieves on linear problems). As it is the case for the other two values of $\bar b$, the rate that may be proved under the assumption $\bar P\succeq 0$ is unduly pessimistic.

By comparing the three panels, we see that the value of the friction parameter that leads to a faster decay
of the ODE solution is $\bar b=2$, i.e.\ the best choice for quadratic objective functions. In this regard, we note that once $t$ is so large that $x(t)$ is close to $x^\star$, the objective function becomes approximately quadratic $f(x) \approx f(x^\star)+(1/2) (x-x^\star)^\T H (x-x^\star)$, with $H$ given by the Hessian matrix of $f$ evaluated at the minimizer.

\subsection{Discrete time analysis}
\label{subsec:nest}
We will now study optimization methods of the form \eqref{eq:nest} for $\alpha_{k}=\alpha$ and $\beta_{k}=\beta$. In order to easily relate what follows to the time-continuous case,  we first introduce as a new variable
 the divided difference, \(k=1, 2,\dots\),
$$%\begin{equation}\label{eq:dd}
d_k = \frac{1}{\delta}(x_k-x_{k-1}),
$$%\end{equation}
where the steplength $\delta = \sqrt{m\alpha}$ is nondimensional ($\beta$ is also nondimensional). \modp{Note this change of variables can be written as
	\begin{equation*}
		\left[\begin{matrix}
			d_k\\
			x_k
		\end{matrix}\right] =\left[\begin{matrix}
			-\frac{1}{\delta} & \frac{1}{\delta}\\
			0 & 1
		\end{matrix}\right]\left[\begin{matrix}
			x_{k-1}\\
			x_k
		\end{matrix}\right] .
	\end{equation*}
}
With the new variable,  \modp{\eqref{eq:nest_cont_framework}} becomes  (\(k=0,1,\dots\))
\begin{subequations}\label{eq:nest3}
\begin{align}
d_{k+1} &=\beta d_k -\frac{\alpha}{\delta} \nabla f(y_k),\\
x_{k+1}  &= x_k+\delta\beta d_k -\alpha \nabla f(y_k) ,\\
y_k &= x_k+\delta\beta d_k,
\end{align}
\end{subequations}
 and these equations are of the form \eqref{eq:control_disc} with $\xi_k= [d_k^{\T},x_k^{\T}]^{\T} \in\R^{2d}$ and
\begin{equation} \label{eq:nest_form}
A = \left[\begin{matrix}\beta I_d & 0\\ \delta\beta I_d& I_d\end{matrix}\right],\quad
B = \left[\begin{matrix}-(\alpha/\delta) I_d\\ -\alpha I_d\end{matrix}\right],\quad
C = \left[\begin{matrix} \delta\beta I_d & I_d\end{matrix}\right],\quad
E = \left[\begin{matrix} 0 & I_d\end{matrix}\right].
\end{equation}

\modg{We now have the following theorem, where, as distinct from the continuous case, the rate of convergence is not given explicitly but as the solution of an algebraic equation.}

\begin{figure}
	\centering
	\begin{subfigure}[b]{0.32\textwidth}
		\centering
		\includegraphics[width=\textwidth]{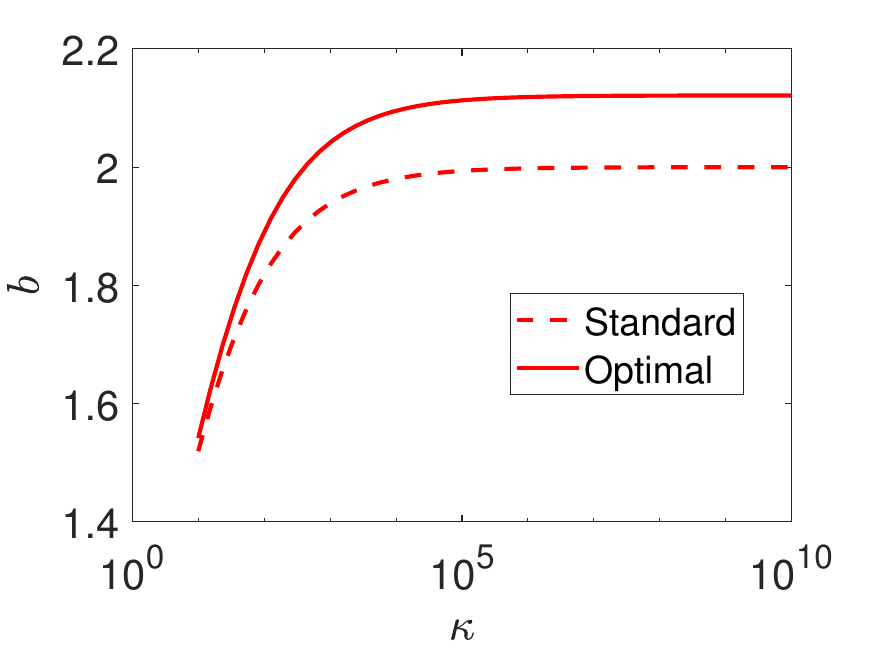}
		\caption{ Parameter choice }
		\label{fig:boptimal}
	\end{subfigure}
	\hfill
	\begin{subfigure}[b]{0.32\textwidth}
		\centering
		\includegraphics[width=\textwidth]{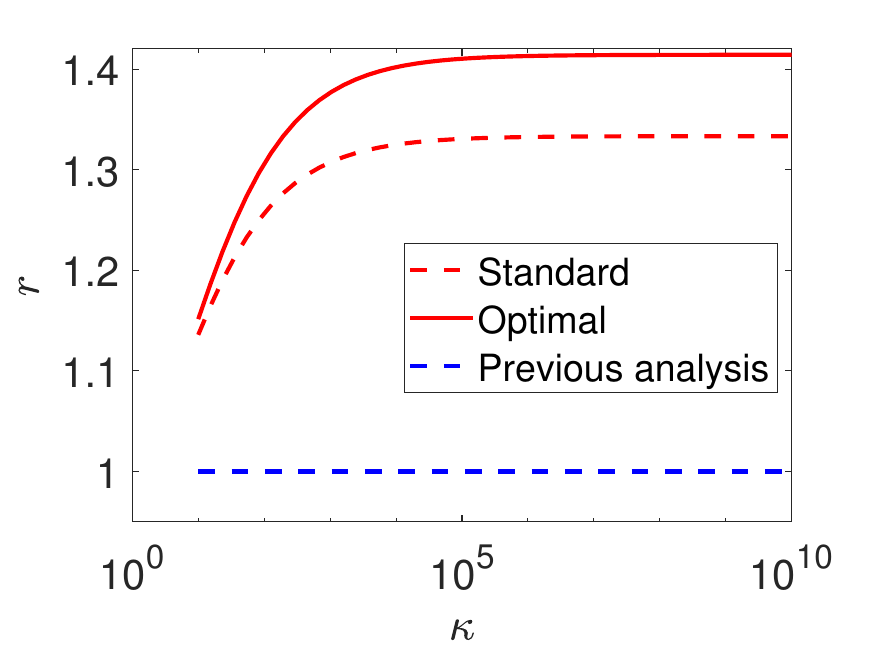}
		\caption{Convergence rate}
		\label{fig:roptimal}
	\end{subfigure}
	\hfill
	\begin{subfigure}[b]{0.32\textwidth}
		\centering
		\includegraphics[width=\textwidth]{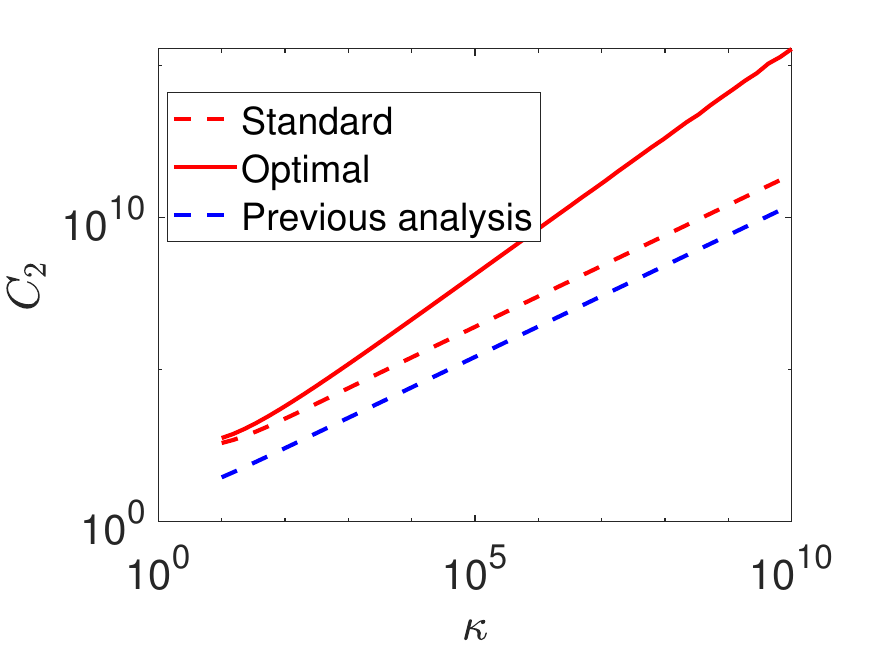}
		\caption{bound $C_2$}
		\label{fig:Coptimal}
	\end{subfigure}
	\caption{\modp{Here we plot the convergence rate obtained by Theorem~\ref{thm:nestconvergence} for the standard choice of $\beta=(\sqrt{\kappa}-1)/(\sqrt{\kappa}+1)$ in the dashed lines and for the optimal choice $\beta=1-b\delta$ derived from the estimate in Theorem~\ref{thm:nestconvergence} in the solid lines. In panel \ref{fig:boptimal} we give the value of $b$ for different values of $\kappa$, in panel \ref{fig:roptimal} we give $r$ where $\rho^2=1-r\delta$, and in panel \ref{fig:Coptimal} we give the constant $C_{2}$ in the estimate
\eqref{eq:comparison_estimate}.}}
\label{fig:optimal_rates}
\end{figure}

\modp{
	\begin{theorem}\label{thm:nestconvergence}
		Fix $0<m<L$, $f\in \mathcal{F}_{m,L}$, and let $\xi_k=[d_k^\T,x_k^\T]^\T$ be the sequence given by \eqref{eq:nest} with parameters $0 \leq \alpha\leq 1/L$, $\beta=1-b\delta$ for  $b>0$, $\delta=\sqrt{m\alpha}$ and initial condition $\xi_{0}=(d^{T}_0,x^{T}_0)$. Then
		\begin{equation} \label{eq:disc_nest_convergence}
			\norm{x_k-x^{\star}}^{2}  \leq C\rho^{2k}\left(f(x_0)-f(x^\star)+\norm{\xi_0-\xi^{\star}}_{\tilde{P}}^2\right).
		\end{equation}
		with
		\begin{eqnarray} \label{eq:p22_disc}
 C&=&\frac{1}{\min{\sigma(\tilde{P})}}, \nonumber \\
		\tilde{P}&=&\frac{m}{2}\left[\begin{matrix}
		p_{22} \delta^{2} -2r\delta+1 & r-\delta p_{22} \\
		r-\delta p_{22}  &p_{22}+1
		\end{matrix}\right] \otimes I_{d},\\
		   p_{22}&=&r\frac{b^2\delta^3 - b^2\delta - 2rb\delta^3 + 2rb\delta + 3r\delta^2 - 2\delta -r}{(2\delta r - 2)},\nonumber
		\end{eqnarray}
		and $\rho^2=1-r\delta$. Here $r$ is the largest positive root of the  polynomial
		\begin{align}  \label{eq:r_disc}
			0&=r\left(1-p_{22}\right)\left(2b+\delta +\delta  p_{22}-3  r+2 \delta   r^2- \delta ^2 p_{22} r+b^2 \delta ^3 -2 b \delta ^2 -b^2 \delta  \right)  \\
			&\qquad\qquad-{\left( p_{22}+ r^2-b  r-\delta   r- \delta  p_{22} r+b \delta ^2  r\right)}^2, \nonumber
%			\varphi(r,b;\delta)&=-b^4 \delta ^6 r+2 b^4 \delta ^4 r-b^4 \delta ^2 r+4 b^3 \delta ^6 r^2-8 b^3 \delta ^4 r^2+4 b^3 \delta ^2 r^2-4 b^2 \delta ^6 r^3-10 b^2 \delta ^5 r^2+8 b^2 \delta ^4 r^3\\
%			&+12 b^2 \delta ^4 r+16 b^2 \delta ^3 r^2-4 b^2 \delta ^3-4 b^2 \delta ^2 r^3-16 b^2 \delta ^2 r-6 b^2 \delta  r^2+4 b^2 \delta +4 b^2 r+12 b \delta ^5 r^3\\
%			&-24 b \delta ^3 r^3-16 b \delta ^3 r+8 b \delta ^2 r^2+8 b \delta ^2+12 b \delta  r^3+16 b \delta  r-8 b r^2-8 b-9 \delta ^4 r^3+8 \delta ^3 r^2\\
%			&+22 \delta ^2 r^3+4 \delta ^2 r-4 \delta  r^4-32 \delta  r^2-4 \delta +3 r^3+12 r			
		\end{align}				
and it must be assumed that the following constraints hold:		
\begin{subequations}	 \label{eq:p22_const}	
\begin{align}
			0&< 1+\frac{p_{22}}{2}-\frac{\sqrt{\delta ^2+1} \sqrt{\delta ^2 {p_{22}}^2-4 \delta   p_{22} r+4r^2+{p_{22}}^2}}{2}+\frac{\delta ^2 p_{22}}{2}-\delta   r,\\
			0&<1+\frac{p_{22}}{2}+\frac{\sqrt{\delta ^2+1} \sqrt{\delta ^2 {p_{22}}^2-4 \delta   p_{22} r+4r^2+{p_{22}}^2}}{2}+\frac{\delta ^2 p_{22}}{2}-\delta   r, \\
			0 & \leq 2b+\delta +\delta  p_{22}-3  r+2 \delta   r^2- \delta ^2 p_{22} r+b^2 \delta ^3 -2 b \delta ^2 -b^2 \delta, \\
			0&\leq 1-p_{22}.
		\end{align}
\end{subequations}				
	\end{theorem}		
}
\modkz{
\begin{proof}[Proof of Theorem \ref{thm:nestconvergence}]
The proof follows from Theorem \ref{theo:main_disc} after setting (see Remarks \ref{rem:choice1a}, \ref{rem:choice1b}  for further explanation)
\begin{equation} \label{eq:choices2a}
a_{0}=1, \qquad \rho^{2}=1-r\delta, \qquad  \beta=1-b\delta, \qquad\ell=0
\end{equation}
and
\begin{equation} \label{eq:choices2b}
P=\frac{m}{2}\left[\begin{matrix}
		p_{22} \delta^{2} -2r\delta+1 & r-\delta p_{22} \\
		r-\delta p_{22}  &p_{22}
		\end{matrix}\right] \otimes I_{d},
\end{equation}
where $p_{22}$ is given by \eqref{eq:p22_disc}. The matrices $A, B, C$ are given in \eqref{eq:nest_form} and then the matrix $\bar{T}$ becomes $T=\widehat{T} \otimes I_{d}$ where
\begin{eqnarray*}
\widehat{t}_{11} &=& -\frac{\delta m}{2}\left(2b+\delta +\delta  p_{22}-3  r+2 \delta   r^2- \delta ^2 p_{22} r+b^2 \delta ^3 -2 b \delta ^2 -b^2 \delta \right), \\
\widehat{t}_{12} &=& \frac{\delta m}{2} \left(p_{22}+ r^2-b  r-\delta   r- \delta  p_{22} r+b \delta ^2  r \right), \\
\widehat{t}_{13} &=& 0, \\
 \widehat{t}_{22} &=& \frac{\delta m}{2} r (p_{22}-1),  \\
 \widehat{t}_{23} &=& 0, \\
 \widehat{t}_{33} &=&  \frac{1}{2}\alpha(L\alpha-1).
\end{eqnarray*}
As distinct to the continuous case,  $\widehat{t}_{33}$ is non-zero so that, in order to ensure that $T \preceq 0$, one needs that
\[
0 \leq \alpha \leq \frac{1}{L}.
\]
In addition, we need to ensure that the matrix $\tilde{P}$ is non-singular, which gives rise to conditions (\ref{eq:p22_const}a) and (\ref{eq:p22_const}b). Finally,  we need the diagonal elements of $\widehat{T}$ to be non-positive thus obtaining conditions (\ref{eq:p22_const}c) and (\ref{eq:p22_const}d), while  the determinant of the leading $2\times 2$ submatrix of $\widehat{T}$ has to be non-negative giving rise to
\begin{align*}
0& \leq r\left(1-p_{22}\right)\left(2b+\delta +\delta  p_{22}-3  r+2 \delta   r^2- \delta ^2 p_{22} r+b^2 \delta ^3 -2 b \delta ^2 -b^2 \delta  \right) \\
			&\qquad\qquad-{\left( p_{22}+ r^2-b  r-\delta   r- \delta  p_{22} r+b \delta ^2  r\right)}^2.
\end{align*}
In the continuous case (which is recovered by setting \(\delta = 0\)), for $b\neq {3\sqrt{2}}/{2}$ the constraint on $\tilde{P}$ is not active  (the active constraint arises from the determinant of the leading submatrix). Since for ill-conditioned problems $\delta \ll 1$, we obtain that $r$ has to correspond to the largest positive root of  \eqref{eq:r_disc} subject to the constraints (\ref{eq:p22_const}a)--(\ref{eq:p22_const}d).
\end{proof}
}
\modg{
\begin{remark} \label{rem:choice1a}
Writing $\rho^2$ and $\beta$ as in terms of new parameters $r$ and $b$ as in \eqref{eq:choices2a} simplifies the algebra and, in addition, prepares the way for the comparison between the discrete and continuous analysis. Similar to what was done in the continuous case, we  set
$\ell=0$, as this does not have  a significant impact on the value of \(\rho\) that results from the analysis. This in turn allows us to simplify things further, since then $T$ is homogeneous in $P$ and $a_{0}$ and we may set $a_{0}=1$.
\end{remark}
}
\modg{
\begin{remark} \label{rem:choice1b}
The rationale behind the choice of the matrix $P$ in \eqref{eq:choices2b} with $p_{22}$ as in \eqref{eq:p22_disc} is similar to the continuous case one. We use the Kronecker product structure and then we  choose the elements of the corresponding $2\times2$ matrix to get $\widehat{t}_{13}=\widehat{t}_{23}=0$, which gives us \eqref{eq:choices2b}. Finally, $p_{22}$ is obtained by setting  $p_{22}$ to be the solution of  $\partial_{p_{22}}(\widehat{t}_{11}\widehat{t}_{22}-\widehat{t}_{12}^{2})=0$. The main difference with the continuous case is that now it is impossible to write explictly \(p_{22}\).
\end{remark}
}
\modkz{
Theorem \ref{thm:nestconvergence} provides a mechanism to determine the convergence rate  $\rho^{2}$ of the algorithm \eqref{eq:nest3} for different choices of $\alpha, \beta$, for  given condition number $\kappa$. One has to solve numerically \eqref{eq:r_disc} for $r$ and check that the constraints (\ref{eq:p22_const}a)--(\ref{eq:p22_const}d) are satisfied.
In particular, in Figure  \ref{fig:fig2} we plot the numerically found value of $r$ as given by Theorem  \ref{thm:nestconvergence} with $\alpha=1/L, \beta=1-{b}/{\sqrt{\kappa}}$ for two choices of $L$ when $m=1$.\footnote{Recall $L$ and $m$ individually may be rescaled. It is the condition number $\kappa=L/m$ that matters.} In addition, we  compare with the analogous curve obtained in \cite{SSKZ21} under the  constraint $ P \succeq 0 $ required by the framework in \cite{FRMP18}. In \cite{SSKZ21}  the best achievable rate is $r = 1$. As we may see, by changing the constraint on $P$  it is possible to prove a significantly better convergence rate.  In particular, for the modified constraint in the present analysis, one can show easily that $b$ may be chosen to get $r=\sqrt{2}-\mathcal{O}(1/\sqrt{\kappa})$, which in turn implies that  in \eqref{eq:conv_disc}:
\[
\rho^{2}=1-\frac{\sqrt{2}}{\sqrt{\kappa}}+\mathcal{O}\left(\frac{1}{\kappa}\right), \quad \kappa \rightarrow \infty.
\]
This is further illustrated in Figure \ref{fig:roptimal} where we plot the value of $r$ as a function of $\kappa$ for the best choice of $\beta$ and $\alpha=1/L$.  Note that the best choice of $\beta$ is given in Figure \ref{fig:boptimal} and it differs from the standard choice \cite{N14}  of $(\sqrt{\kappa}-1)/(\sqrt{\kappa}+1)$
plotted using the dashed line.
}
\modkz{
In order to make comparisons between different Lyapunov functions we convert our estimates to bound $\norm{x_{k}-x^{\star}}^{2}$ in terms of $\norm{\xi_{0}-\xi_{*}}^{2}$, i.e.
\begin{equation}\label{eq:comparison_estimate}
	\norm{x_k-x^\star}^2 \leq C_2 \rho^{2k} (\norm{x_0-x^{\star}}^2+\norm{d_0}^2).
\end{equation}
We plot $C_{2}$ as a function of $\kappa$ in Figure \ref{fig:Coptimal}. As we can see for the standard choice of $\beta$ $C_{2}$ grows in a similar manner as obtained in the analysis in \cite{SSKZ21}. On the other hand, if one uses the  $\beta$ that maximizes $\rho^{2}$, then the constant $C_{2}$ grows faster. This is  to be expected since $\tilde{P}$ becomes singular as $r$ approaches the optimal value.
}

\section{Connecting optimization algorithms and Polyak's ODE}
\label{sec:opt}

We now discuss the relations between the continuous and discrete time studies presented above.

\subsection{The Nesterov algorithm as an integrator}
\label{subsec:opt1}
For suitable parameter choices, the Nesterov algorithm \eqref{eq:nest} is a discretization of Polyak differential equation. However, as discussed in detail in \cite{SSKZ21}, such a discretization does not correspond to any of the more familiar classes of ODE solvers, such as linear multistep or Runge Kutta (RK) methods. In particular we remark that in \eqref{eq:nest} $\nabla f$ is not evaluated at the approximations $x_k$ delivered by the algorithm.
As we shall see presently, it turns out that the Nesterov algorithm is an example of the class of
Additive Runge-Kutta (ARK) algorithms, a generalization of the RK integrators considered by several authors after its introduction by Cooper \cite{Cooper80,Cooper83}.

Additive Runge-Kutta (ARK) algorithms integrate systems of differential equations $(d/dt) z = g(z)$ in cases where it makes sense to decompose $g(z)$ as a sum $g(z) = \sum_{\nu = 1}^Ng^{[\nu]}(z)$. In the plain RK case, the numerical solution is advanced over a time step $z_k\mapsto z_{k+1}$ by evaluating $g(z)$ at a sequence of so-called stage vectors $Z_{k,1}$, \dots, $Z_{k,s}$ and then setting $z_{k+1} = z_k +\sum_{i=1}^s b_i g(Z_{k,i})$, where the $b_i$ are suitable weights. In turn, for the explicit algorithms we are interested in, the stages are computed successively,  $i = 1, \dots, s$, as $Z_{k,i} = z_k + h\sum_{j=1}^{i-1} a_{i,j}g(Z_{k,j})$, with suitable coefficients $a_{i,j}$. ARK algorithms are entirely similar, but evaluate the individual pieces $g^{[\nu]}(z)$ rather than $g(z)$.

With $z=[v^\T,x^\T]^\T\in\R^{2d}$, the system \eqref{eq:ode1} may be rewritten as
$$
\frac{d}{dt} z= g^{[1]}(z) + g^{[2]}(z) + g^{[3]}(z):=
\begin{bmatrix} -\bar{b} \sqrt{m} v\\0 \end{bmatrix}+
\begin{bmatrix} -\frac{1}{\sqrt{m}}\nabla f(x)\\0 \end{bmatrix}+
\begin{bmatrix} 0\\\sqrt{m} v\end{bmatrix};
$$
the three parts of $g(z)$ respectively represent the friction force, potential force and inertia in the oscillator. It is easily checked that, if we choose a steplength $h>0$, and see $d_k$ and $x_k$ as approximations to $v(kh)$ and $x(kh)$ respectively, then a step $(d_k,x_k)\mapsto (d_{k+1},x_{k+1})$ of the optimization algorithm \eqref{eq:nest3}
with parameters $\alpha = h^2$, $\beta = 1-h\bar{b}\sqrt{m}$, $\delta = \sqrt{m} h$ is just one step
$z_k\mapsto z_{k+1}$ of the ARK integrator for \eqref{eq:ode1} given by:
\begin{eqnarray*}
Z_{k,1} &=& z_k,\\
Z_{k,2} & =& z_k+ h g^{[1]}(Z_{k,1}),\\
Z_{k,3} & =& z_k+h  g^{[1]}(Z_{k,1})+ h g^{[3]}(Z_{k,2}),\\
Z_{k,4} & =& z_k+h  g^{[1]}(Z_{k,1})+ h g^{[3]}(Z_{k,2})+ h g^{[2]}(Z_{k,3}),\\
z_{k+1} &=& z_k + h g^{[1]}(Z_{k,1})+ h g^{[2]}(Z_{k,3})+ h g^{[3]}(Z_{k,4}).
\end{eqnarray*}
The stage vectors have $Z_{k,1} = [d_k^\T,x_k^\T]^\T$, $Z_{k,2} = [\beta d_k^\T,x_k^\T]^\T$,
$Z_{k,3} = [\beta d_k^\T,y_k^\T]^\T$, $Z_{k,4} = [d_{k+1}^\T,y_k^\T]^\T$, and therefore
the computation of the second, third and fourth stages incorporate successively the contributions of friction, inertia and potential force.

If we now think that the value of $h>0$ varies and consider the optimization algorithm \eqref{eq:nest3} with $\alpha = h^2$, $\beta = 1-h\bar{b}\sqrt{m}$, the standard theory of numerical integration of ODEs shows that, if the initial points $x_{-1}$ and $x_0$ are chosen in such a way that, as $h\rightarrow 0$, $x_0$ and $(1/h)(x_0-x_{-1})$ converge to limits $A$ and $B$, then, in the limit of $kh\rightarrow t$, $x_k$ and $(1/h)(x_{k+1}-x_k)$ converge to $x(t)$ and $\dot x(t)$ respectively, where $x(t)$ is the solution of \eqref{eq:polyak} with initial conditions $x(0) =A$ and $\dot x(0)=B$. In addition, the discrete Lyapunov function of the optimization algorithm in Section 3 may be easily seen to converge \modp{pointwise} to the Lyapunov function of the ODE found in this section, \modp{see \cite{SSKZ21}}. Finally the discrete decay factor over $k$ steps  $(1-\sqrt{m} r h)^k$ converges to the continuous decay factor $\exp(-\lambda t)$. These facts in particular explain that, in Figure~\ref{fig:conta}, the graph of the relation between $\bar b$ and $\bar r$ that holds for the ODE is indistinguishable from the corresponding graph for the optimization algorithm  when $\kappa$ is large ($\kappa$ being large corresponds to $h$ being small).

\subsection{Discretizations that do not succeed in getting acceleration}
Many recent contributions have derived optimization algorithms by discretizing  suitably chosen dissipative ODEs. It is well known that, unfortunately, many properties of ODEs are likely to be lost in the discretization process, even if high-order, sophisticated integrators are used. The archetypical example is provided by the discretization of the standard harmonic oscillator: most numerical methods, regardless of their accuracy, provide solutions that either decay to the origin or spiral out to infinity as the number of computed points grows unboundedly. Similarly, discretizations of \eqref{eq:polyak} are likely not to share the favourable decay properties in Section \ref{subsec:polyak_ODE}.

Let us consider the following extension of the optimization algorithm \eqref{eq:nest}:
\begin{subequations}\label{eq:nest4}
\begin{align}
x_{k+1} &=y_k+\beta(x_k-x_{k-1})-\alpha\nabla f(y_k),\\
y_k  &= x_k+\gamma  (x_k-x_{k-1}),
\end{align}
\end{subequations}
with the additional parameter $\gamma$. The choice $\gamma = 0$ yields the \emph{heavy ball} algorithm, which (see \cite{SSKZ21}) corresponds to a \lq\lq natural\rq\rq\ standard linear multistep discretization of the Polyak equation \eqref{eq:polyak} where $\nabla f$ is evaluated at the approximations $x_k$. Unfortunately the heavy ball algorithm  does not provide acceleration. As shown in \cite{SSKZ21} for  $\gamma=0$ (or more generally for $\gamma\neq \beta$), the optimization algorithm \eqref{eq:nest4} does not inherit a Lyapunov functions from the Polyak ODE. The analysis in that paper hinges on a study of the nondimensional quantity $c : = t_{11}/(m\delta)$, which for $\widehat{T}\preceq 0$ has to be $\leq 0$ and for a discretization of an ODE has a finite limit as $\delta\rightarrow 0$.
When $\gamma=0$,  the expression for the quantity $c$ includes a positive contribution $\delta (\kappa-1)\beta^2/2$; for acceleration, $\delta$ has to be $\mathcal{O}(1/\sqrt{\kappa})$ which makes it impossible for $c$ to be $\leq 0$.

The unwelcome presence of $\kappa$ in $t_{11}$ may be traced back to the appearance of $L$ in the matrix $N^{(1)}$ in Theorem~\ref{theo:main_disc}. Nesterov's algorithms of the family \eqref{eq:nest} do not suffer from that appearance because for them the matrix $EA-C$ that multiplies $(L/2)I_d$ in the recipe for $N^{(1)}$ vanishes. The condition $EA-C=0$ appears then to be of key importance in the success of Nesterov algorithms; we put it into words by saying that one has to impose that the point $y_k =C\xi_k$ where the gradient is evaluated  has to coincide with the point $x_{k+1} = EA \xi_k$ that the algorithm would yield if $u_k = \nabla f(y_k)$ happened to vanish (see  \eqref{eq:control_disc}). This suggests that the integrator has to treat the potential force and the friction force in the oscillator separately, something that may be achieved by ARK algorithms, but not by more conventional linear multistep or RK methods that do not avail themselves of the separate pieces $g^{[1]}(z)$, $g^{[2]}(z)$, $g^{[3]}(z)$ but are rather formulated in terms of $g(z)$.

\section{Derivation and analysis of a new second order ODE}
\label{sec:polyakplus}

In the last few years  there has been a number of works that have studied ways of accelerating convergence towards equilibrium for dynamics of   stochastic differential equations \cite{LNP13,DNP17,FTN21a,HHS05,HNW15} . When the dynamics of the underlying SDE are linear this problem is directly connected to finding the minimum of a quadratic function $f(x)=(1/2)x^{\T}Sx+c^{\T}x$, with $S=S^\T \succ 0$. For simplicity we will assume that $c=0$ and  in this case the GF \eqref{eq:gf} obtains the simple form
\modp{\[
\frac{dx(t)}{dt}=-Sx(t)
\]}
and the speed of convergence towards zero is dictated by the minimum eigenvalue of $S$. It is possible to increase the speed of convergence towards zero by introducing a \emph{non-reversible} perturbation to the above equation. More precisely, it is easy to show  \cite{LNP13} that the dynamics of
\modp{
	\[
\frac{dx(t)}{dt}=-(I+J)Sx(t), \qquad J=-J^{\T}
\]
}
 yields  faster  convergence towards zero than the original GF. Furthermore, as discussed in \cite{LNP13} there is an optimal perturbation $J^{\star}$ for which the rate of convergence towards zero is maximized, with the maximum value being { $\text{Tr}(S)/d$} where $d$ is the dimension of the matrix. A natural question to ask is if this kind of acceleration remains true when  $f \in \mathcal{F}_{m,L}$ is not quadratic.
 In this case the perturbed ODE has the form
 \modp{	\[
 	\frac{dx(t)}{dt}=-(I+J)\nabla f(x(t)), \quad J=-J^{\T}.
 	\]}
 	This equation and discretizations of it were studied in \cite{FTN21b}. In particular, it was shown that upon assuming additional information about the eigenvalues of the Hessian of $f$, convergence rates may improve both in the continuous and discrete setting. Here we will instead consider the GD \eqref{eq:gf} for an appropriately chosen extended objective function and modify  its dynamics with a simple non-reversible perturbation.
 %
% It turns out that,  working directly with the objective function $f(x)$, it does not appear to be easy to \modk{quantify the improved convergence rate towards} equilibrium of the perturbed ODE
%
%\[
%\frac{dx}{dt}=-(I+J)\nabla f(x), \quad J=-J^{\T},
%\]
%when $f \in \mathcal{F}_{m,L}$ without additional  hypothesis on $f$ \cite{FTN21b}.
%However, if we consider
\modk{In this case it is possible to fully quantify the increase in the convergence rate without any additional assumptions on $f$}.

\modk{We introduce} an auxiliary variable $y\in \R^d$ and the extended objective function $F(y,x) = (L/2)\norm{y}^2+f(x)$ with minimum at $(0,x^\star)$. The corresponding GF is
\modp{$$
\frac{d}{d\tau}\left[\begin{matrix}y(\tau)\\ x(\tau) \end{matrix}\right] = - \left[\begin{matrix} Ly(\tau)\\\nabla f(x(\tau))\end{matrix}\right]
$$}
and we perturb the right hand-side by adding a skew-symmetric term to get
\modp{$$
\frac{d}{d\tau}\left[\begin{matrix}y(\tau)\\ x(\tau) \end{matrix}\right] =
 - \left[\begin{matrix}Ly(\tau)\\\nabla f(x(\tau))\end{matrix}\right]
+ K \left[\begin{matrix}0&-I\\I&0\end{matrix}\right] \left[ \begin{matrix}Ly(\tau)\\\nabla f(x(\tau))\end{matrix}\right]
$$}
where $K\geq 0$ is a perturbation parameter. For $K=0$, $x$ evolves as in \eqref{eq:gf} (but the time variable here has been relabelled for reasons that will become clear immediately).
By replacing the variables $y$ and $\tau$ and the parameter $K$ by $v$, $t$ and $\bar b\geq 0$ respectively, with
$$
   y = \sqrt{\frac{m}{L}} v, \qquad
   \tau = \frac{\bar b \sqrt{m}}{L} t, \qquad
   K = \frac{1}{\bar b}  \sqrt{\frac{L}{m}},
$$
the system becomes
\begin{subequations}\label{eq:polyakplus}
\begin{align}
\frac{d}{dt} v\modp{(t)} &=  -\bar b \sqrt{m} v\modp{(t)}-\frac{1}{\sqrt{m}} \nabla f(x\modp{(t)}),\\
\frac{d}{dt} x\modp{(t)}&= \frac{\bar b\sqrt{m}}{L} \nabla f(x\modp{(t)})+\sqrt{m} v\modp{(t)}.
\end{align}
\end{subequations}
Comparing these expressions with \eqref{eq:ode1}, we see that we are dealing here with a perturbation of Polyak's equation, where now $\nabla f(x)$ is used both in the $v$ and $x$ equations; Polyak equation is retrieved in the limit $L\uparrow \infty$ with fixed $m$. For this reason we shall refer to \eqref{eq:polyakplus} as the \emph{Polyak+} system. As noted before, as the friction coefficient $\bar{b}$ grows unboundedly (i.e.\ $K\downarrow 0$) with fixed $L$ and $m$, the dynamics of $x$ under \eqref{eq:polyakplus} approaches GD; on the other hand, in the limit $\bar{b} \downarrow 0$, \eqref{eq:polyakplus} becomes a Hamiltonian (nondissipative) system.

 The system \eqref{eq:polyakplus} is easily cast in the control framework of Section~\ref{sec:pre} and we use Theorem~\ref{theo:main_con} to investigate to what extent it improves on Polyak's ODE.
The $\bar{t}_{ij}$ in \eqref{eq:T} are found to be
\begin{align*}
  \bar t_{11} &=  -2\bar b \sqrt{m} \bar p_{11}+2\sqrt{m} \bar p_{12}+\lambda \bar p_{11},\\%
  \bar t_{12} &= -\bar b \sqrt{m} \bar p_{12}+\sqrt{m} \bar p_{22}+\lambda \bar p_{12}, \\%
  \bar t_{13} &= -(1/\sqrt{m})\bar p_{11}+\sqrt{m}/2-\bar{b} \sqrt{m}\bar{p}_{12}/L, \\%
  \bar t_{22} &= \lambda  \bar p_{22}-(m/2)\lambda,\\%
  \bar t_{23} &= -(1/\sqrt{m})\bar p_{12}+\lambda/2-\bar{b}\sqrt{m}\bar p_{22}/L,\\%
  \bar t_{33}&=  - \bar{b}\sqrt{m}/L.
\end{align*}
To carry out the analysis it is convenient to introduce $\zeta = 1/L$; the limit value $\zeta= 0$ then corresponds to Polyak's ODE. We saw in \modkz{Theorem~\ref{thm:polyakconvergence}} that, in the $\zeta = 0$ case, each rate $\bar r<\sqrt{2}$ may be achieved with two different values of $\bar b$, one below $3\sqrt{2}/2$ and the other above. We  \modkz{now perform the analysis of \eqref{eq:polyakplus}} for $\bar b < 3\sqrt{2}/2$. When determining $\lambda=\sqrt{m} \bar r$ and the elements $\bar{p}_{ij}$ we  \modkz{impose that
\begin{equation} \label{eq:step1}
\bar p_{11} = m/2
\end{equation}
and  operate under the assumption that the matrix $\widehat{\bar{T}}$ has rank $\leq 1$.}
%\item We have\begin{equation}\label{eq:step1}\bar p_{11} = m/2.\end{equation}
%\end{enumerate}
In extensive experimentation we have observed that these two conditions hold when $\lambda$ is numerically maximized subject to the constraints $\widetilde{P} \succ 0$, $\widehat{\bar {T}} \preceq 0$. \modkz{Note also that they are satisfied in Polyak's, $\zeta = 0$, case (see the proof of Theorem \ref{thm:polyakconvergence}).}
%: for the first assumption recall that we saw in Section~\ref{subsec:polyak_ODE} that, for $\bar b < 3\sqrt{2}/2$ , when $\bar r$ is maximized all elements of $\widehat{\bar {T}}$ vanish, except perhaps $\bar t_{22}$ and for the second assumption see \eqref{eq:p11}.
\modkz{We point out that  our assumption  on the rank is equivalent to the requirement that all $2\times 2$ submatrices of $\widehat{\bar{T}}$ are singular}.

The assumptions above  uniquely determine $\lambda$ and $\bar{P}$ in Theorem~\ref{theo:main_con}. We proceed as follows:
\begin{itemize}\item By imposing that the determinant of the first and third rows and columns of $\widehat{\bar {T}}$ vanish we find
\begin{equation}\label{eq:step2}
\bar r = 2\bar{b} - (4/m)\bar{p}_{12}-(2\bar{b}/m)\bar{p}_{12}^2\zeta.
\end{equation}
\item By annihilating the determinant of the second and third rows and first and third columns
\begin{equation}\label{eq:step3}
\bar{p}_{22} = \bar{p}_{12}^2/m.
\end{equation}
\item We take the expressions for $\bar r$ and $\bar{p}_{22}$ just found to the equation $\bar{t}_{22}\bar{t}_{33}-\bar{t}_{23}^2 = 0$. This yields an algebraic relation between $\bar{p}_{12}$ and $\zeta$:
\modkz{\begin{eqnarray}
      \label{eq:step4}
    (2\bar b^2 \bar p_{12}^4+m^2 \bar b^2 \bar p_{12}^2) \zeta^2 &+& (8\bar b \bar p_{12}^3-2m\bar b^2 \bar p_{12}^2+2m^2\bar b \bar p_{12}-m^3\bar b^2 ) \zeta  \\
   &+&  (- 3 \bar p_{12} + m\bar b)^2=0. \nonumber
\end{eqnarray}}
\end{itemize}
 The conditions \eqref{eq:step2}--\eqref{eq:step4} guarantee that the rank of $\widehat{\bar T}$ is $\leq 1$, i.e.\ the matrix has at least two zero eigenvalues. Since $\bar t_{33}<0$ for $\bar b>0$ and $\zeta >0$ the matrix will have rank exactly one and be negative semidefinite.

The algebraic curve \eqref{eq:step4} in the $(\bar p_{12},\zeta)$ plane contains the points $P_1 =(0, 1/m)$ and $P_2 = (m\bar b/3,0)$. The first corresponds to the $L=m$ (i.e.\ $\kappa =1$) situation; the second was known to us at it corresponds to Polyak's equation. The global behavior of the curve \eqref{eq:step4} may be investigated by solving the quadratic equation for $\zeta$. Restricting the attention to $0\leq \bar p_{12} \leq m\bar b/3$, there is a branch of the curve where to each value of $\bar p_{12}$ there corresponds a unique value of $\zeta$ (i.e.\ of $L$), so that as $\bar p_{12}$ increases monotonically from $0$ to $m\bar b/3$, $\zeta$ decreases monotonically from $1/m$ to $0$ (or $\kappa$ increases from $1$ to $\infty$). Note that this branch connects the points $P_1$ and $P_2$. Once $\zeta=\zeta(\bar p_{12})$, for given $m$, $\bar b$, has been determined in this way, the relations \eqref{eq:step1}--\eqref{eq:step3} determine $\bar r$, $\bar p_{11}$, $\bar p_{22}$ as functions of $\bar p_{12}\in[0,m\bar b/3]$ with known expressions (that we will not reproduce here). It is easily checked that, for the values of $\bar p_{ij}$ defined in this way, $\widetilde P$ is in fact $\succ 0$. Numerical experiments confirm that maximizing $\lambda$ subject to the constraints $\widetilde{P} \succ 0$, $\widehat{\bar {T}} \preceq 0$ for different choices of $L$, $m$ and $\bar b$ leads to the values of $\bar p_{ij}$ and $\bar r$ we have just constructed analytically. This confirms that the procedure we have followed succeeds in identifying the best $\lambda$ and $\bar P$ in Theorem~\ref{theo:main_con} or, in other words, that the two assumptions  formulated at the outset are valid.

\begin{figure}[h]
     \centering
\includegraphics[scale=0.5]{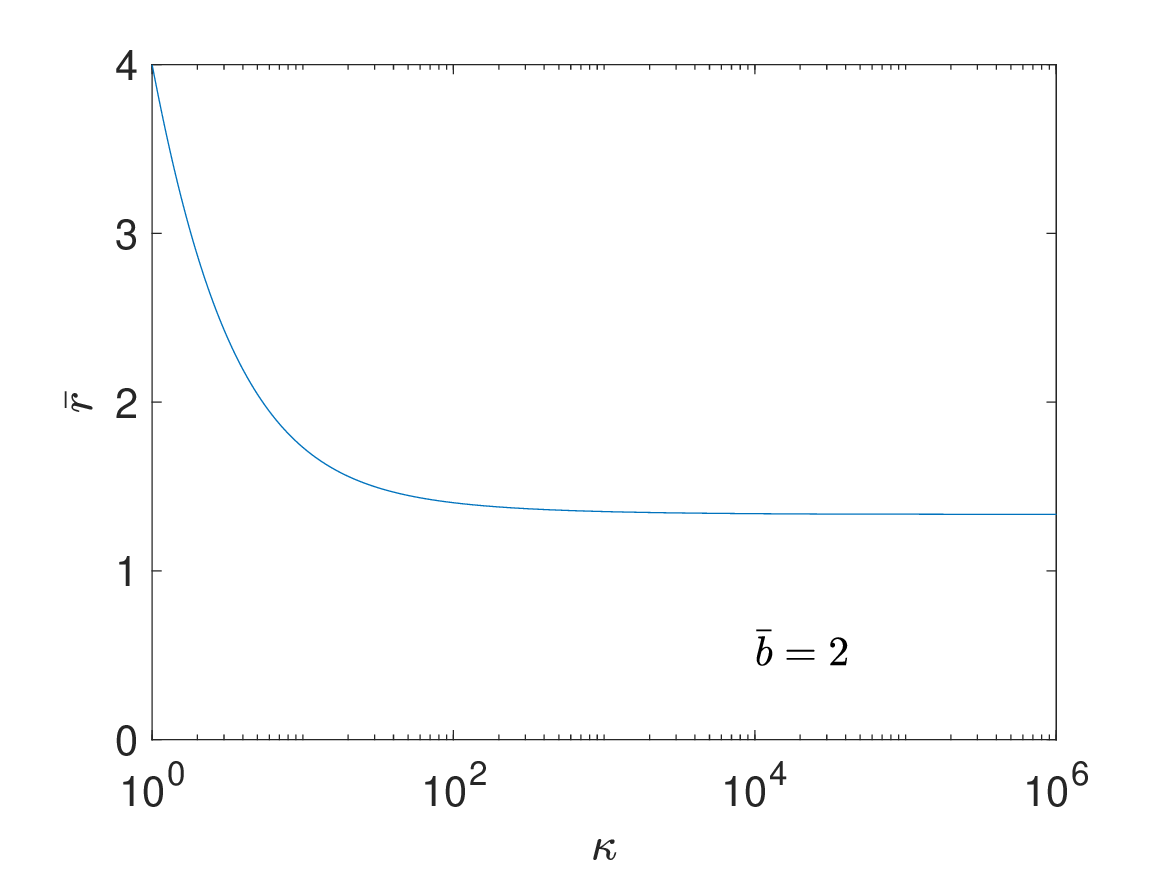}
\caption{Converge rate $r$ as a function of the condition number $\kappa$ for  \eqref{eq:polyakplus}}
\label{fig:polyak_plus_rate}
     \end{figure}

In Fig.~\ref{fig:polyak_plus_rate} we have plotted in the $(\kappa,\bar r)$ plane the parametric curve
$\kappa = L/m = 1/(m\zeta(\bar p_{12}))$, $\bar r = \bar r (p_{12})$, $\bar p_{12} \in [0, m\bar b/3]$ when $\bar b = 2$. Although the system \eqref{eq:polyakplus} clearly improves on Polyak's dynamics for small $\kappa$, the improvement becomes negligible as $\kappa$ increases, i.e.\ in the regime where it would really be needed. We now make this matter more precise.
At the point $P_1$, where $\kappa = 1$, \eqref{eq:step2} yields $\bar r = 2\bar b$; this is to be compared with the value
$\bar r =2\bar b/3$ obtained in Section~\ref{sec:polyak} for Polyak's equation. Thus, the introduction of gradient in the first equation in \eqref{eq:polyakplus} increases $\lambda$ by a factor of $3$.

Let us now consider the neighbourhood of the point $P_2$, i.e.\ the $\kappa\gg 1$ regime.
By implicit differentiation of \eqref{eq:step4}, we find that, at this point, the Taylor expansion of $\zeta$ as a function of $p_{12}$ is given by
$$
\zeta = \frac{243}{9 \bar{b}^2-2\bar{b}^4}\big(p_{12}-\frac{\bar b}{3}\big)^2+\mathcal{O}\Big(\big(p_{12}-\frac{\bar b}{3}\big)^3\Big), \qquad \bar{p}_{12} \rightarrow m\bar{b}/3.
$$
On the other hand, using the expression of $\bar r$ as a function of $\bar{p}_{12}$, one finds
$$
\bar{r} = \frac{2\bar{b}}{3} -\frac{4}{m}\left( \bar{p}_{12}-\frac{m\bar{b}}{3}\right)
+\mathcal{O}\left( \bar{p}_{12}-\frac{m\bar{b}}{3}\right)^2,\qquad \bar{p}_{12} \rightarrow m\bar{b}/3
$$
and, combining the last equations, we obtain
after eliminating $\bar{p}_{12}$,
\begin{equation}\label{eq:polyakplusrate}
\bar{r} = \frac{2\bar{b}}{3}+ C(\bar b) \frac{1}{\sqrt{\kappa}}+\mathcal{O}\left(\frac{1}{\kappa}\right),
\qquad \kappa\rightarrow \infty,
\end{equation}
with
$$
C(\bar b) = 4\sqrt{\frac{9\bar{b}^2-2\bar{b}^4}{243}},\qquad
 \bar{b} < 3\sqrt{2}/2.
$$
 Since $C(\bar b) >0$, we conclude that,  for large condition number $\kappa$, the Polyak+ system in fact achieves a rate larger than the rate ${2\bar{b}}/{3}$ for Polyak's ODE. Unfortunately, in order to maximize the leading  term, $2\bar{b}/3$, in the expansion \eqref{eq:polyakplusrate}, $\bar b$ has to be chosen close to the upper limit $3\sqrt{2}/2$ and, as $\bar{b}\uparrow 3\sqrt{2}/2$,   the increment $C(\bar b)/\sqrt{\kappa}$ vanishes. For instance, for $\bar{b} = 2.1$, \eqref{eq:polyakplusrate} becomes $\bar{r} \approx 1.4000+0.2286/\sqrt{\kappa}$; for $\kappa = 10^4$ the increment is only 0.0022.
Therefore, as in the particular case depicted in Fig.~\ref{fig:polyak_plus_rate}, the improvement in rate of the Polyak+ system on the Polyak ODE is indeed negligible, except in the uninteresting case of small $\kappa$. For this reason we have not undertaken the analysis of optimization algorithms based on discretizations of the Polyak+ system.

\section{Stochastic problems: The case of over-parameterized models}
\label{sec:over}
In this section, we extend the Lyapunov function approach to analyse the performance of optimization methods applied to specific modern machine learning models. In particular, we study models such as non-parametric regression or overparameterised deep neural models that are expressive enough to fit or \emph{interpolate} the data set completely \cite{ZBH17, MBB18}. For these models the function $f(x)$ that one is interested in minimising has the following structure
\begin{equation} \label{eq:struct}
f(x)=\frac{1}{N} \sum_{i=1}^{N}f_{i}(x).
\end{equation}
%In the context of supervised learning each of the terms $f_{i}(x)$ corresponds to the
%loss function for the point $(\chi_{i}, \psi_{i})$ when the model parameters are equal
%to $x$. Here $\chi_{i},\psi_{i}$ refer to the feature vector and label for point $i$
%respectively. Some common choices of the loss function include the squared loss
%where $f_{i}(x)=\frac{1}{2}(x^{\T}\chi_{i}-\psi_{i})^{2}$, the hinge loss where $f_{i}
%(x)=\max(0,1-\psi_{i}x^{\T}\chi_{i})$, or the square hinge loss where $f_{i}
%(x)=\max(0,1-\psi_{i}x^{\T}\chi_{i})^{2}$.
Due to the structure of $f$ in \eqref{eq:struct} any gradient based algorithm would need to calculate
\[
\nabla f(x)=\frac{1}{N}\sum_{i=1}^{N} \nabla f_{i}(x)
\]
which when $N$ is large may be computationally very expensive. A typical strategy followed in stochastic optimization algorithms is to replace the gradient with a random unbiased estimator of it. In the simplest possible case, one uses the following estimator
\[
\widehat{\nabla} f(x)=\nabla f_{i_{\omega}},
\]
where $i_\omega$ is a uniform random variable in the set of integers $\{1,\cdots,n\}$. More generally, and without necessarily assuming the finite sum-structure one replaces the full gradient by
\[
\widehat{\nabla} f(x)=\nabla f(x,z),
\]
where $z$ can be thought of as the random gradient noise, which we assume satisfies $\mathbb{E}(\nabla f(x,z))=\nabla f(x)$.

\subsection{A framework for stochastic algorithms}
We consider optimization algorithms with random noise analogously to \eqref{eq:control_disc} with the formulation
\begin{subequations} \label{eq:control_disc_stochastic}
	\begin{align}
		\xi_{k+1} &=A\xi_{k}+B\tilde{u}_{k}, \\
		\tilde{u}_{k} & = \widehat{\nabla} f(y_{k}), \\
		y_{k} & =C\xi_{k}, \\
		x_{k} &=E \xi_{k},
	\end{align}
\end{subequations}
where $\xi_{k} \in \R^{n}$ is the state, $\tilde{u}_{k} \in \R^{d}$ is the random input
$(d \leq n)$, $y_{k}  \in \R^{d}$ is the feedback output that is mapped  to
$\tilde{u}_{k}$ by the random nonlinear map $\widehat{\nabla} f$. We assume here that at each step the random gradient is chosen to be independent of the current state, i.e. $\widehat{\nabla} f(y_{k}) = \nabla f(y_k,z_k)$ for some random variable $z_k$ independent of $y_k$.

\begin{theorem} \label{theo:main_disc_stochastic}
	Suppose that, for \eqref{eq:control_disc_stochastic}, there exist $a_{0} >0$,  $\rho \in (0,1)$,
	and a symmetric matrix $P$, with
	$\widetilde{P}:=P+(a_{0}m/2)E^{\T}E \succ 0$,
	such that
	\begin{equation} \label{eq:LMI_disc_stochastic}
		T= M^{(0)}
		+a_{0}\rho^{2}M^{(1)}+a_{0}(1-\rho^{2})M^{(2)} \preceq 0,
	\end{equation}
	where
	\[
	M^{(0)}=\left[ \begin{matrix}
		A^{\T}PA-\rho^{2}P & A^\T PB \\
		B^{\T}PA & \rho_0 B^{\T} P B
	\end{matrix}\right],
	\]
	and
	\[
	M^{(1)}=N^{(1)}+N^{(2)}, \quad M^{(2)}=N^{(1)}+N^{(3)},
	\]
	with
	\begin{align*}
		N^{(1)} &= \left[\begin{matrix}\frac{L}{2}(EA-C)^\T(EA-C) & \frac{1}{2}(EA-C)^\T(LEB+1) \\ \frac{1}{2}(LEB+1)^\T(EA-C) & \frac{L\rho_0}{2}(EB)^\T EB+\frac{1}{2}(EB+(EB)^\T)\end{matrix}\right],
		\\
		N^{(2)} &= \left[ \begin{matrix}
			C-E & 0 \\
			0 & I_{d}
		\end{matrix}\right]^{\T}\left[ \begin{matrix}
			-\frac{m}{2}I_{d} & \frac{1}{2}I_{d} \\
			\frac{1}{2}I_{d} & 0
		\end{matrix}\right]\left[ \begin{matrix}
			C-E & 0 \\
			0 & I_{d}
		\end{matrix}\right],
		\\
		N^{(3)} &=\left[ \begin{matrix}
			C^{\T} & 0 \\
			0 & I_{d}
		\end{matrix}\right]\left[ \begin{matrix}
			-\frac{m}{2}I_{d} & \frac{1}{2}I_{d} \\
			\frac{1}{2}I_{d} & 0
		\end{matrix}\right]\left[ \begin{matrix}
			C & 0 \\
			0 & I_{d}
		\end{matrix}\right].
	\end{align*}
	Assume there exists $\rho_0>0$ such that for all $y\in \R^d$
	\begin{subequations}\label{eq:SGC}
		\begin{align}
			&\mathbb{E}[\widehat{\nabla} f(y)^\T(EB)^\T(EB)\widehat{\nabla} f(y)]\leq \rho_0 \nabla f(y)^\T(EB)^\T(EB)\nabla f(y),\\
			&\mathbb{E}[\widehat{\nabla} f(y)^\T B^\T PB\widehat{\nabla} f(y)]\leq \rho_0 \nabla f(y)^\T B^\T PB\nabla f(y).
		\end{align}
	\end{subequations}
	Then, for $f \in \mathcal{F}_{m,L}$, $\rho_0> 1$ and $V$ given by \eqref{eq:liap_disc}, the sequence $\{x_{k}\}$ satisfies
	\begin{equation*} %\label{eq:conv_disc_stoch}
		\mathbb{E}[\norm{x_{k}-x^{\star}}^{2}] \leq \max \sigma(E^{\T}E)\mathbb{E} [\norm{\xi_k-\xi^{\star}}] \leq \frac{ \max \sigma(E^{\T}E)}{\min \sigma(\widetilde{P})} V_{0}(\xi_0)\rho^{2 k}.
	\end{equation*}
\end{theorem}
\modp{The proof of this theorem can be found in Section \ref{subsec:proof_main_disc_stochastic}.}

{
\begin{remark}
Conditions \eqref{eq:SGC} are generalisations of the strong growth condition in \cite{VBS19} which is satisfied if there exists $\rho_0>0$ such that
\begin{equation}\label{eq:SGCf}
	\mathbb{E}[\lVert \widehat{\nabla} f(y)\rVert^2]\leq \rho_0\lVert \nabla f(y)\rVert^2.
\end{equation}
Such a condition implies that $\widehat{\nabla}f(x_\ast)=0$ almost surely. \modkz{An example of such function is the square-hinge loss  function when the data is linearly separable (see discussion in \cite{VBS19}).}
\end{remark}
}

\subsection{A family of stochastic optimization algorithms}
{We now consider the following family of stochastic optimization algorithms}
\begin{subequations}\label{eq:Vaswani_sys}
	\begin{align}
		x_{k+1}&=y_k-\eta \widehat{\nabla} f(y_k),\\
		y_k&= \tilde{\alpha} v_k+(1- \tilde{\alpha})x_k,\\
		v_{k+1}&= \tilde{\beta} v_k+(1- \tilde{\beta})\zeta_k-\gamma \eta \widehat{\nabla} f(y_k).
	\end{align}	
\end{subequations}
This family was considered in \cite{VBS19} as a generalisation of the accelerated coordinate descent method \cite{N12}. By introducing the variable $d_k=v_k-w_k$ we can write the system \eqref{eq:Vaswani_sys} in a form similar to \eqref{eq:nest3} as follows:
\begin{subequations}\label{eq:Nesterov_stoch_sys}
	\begin{align}
		d_{k+1}&=(1- \tilde{\alpha}) \tilde{\beta} d_k-\eta (\gamma-1)\widehat{\nabla}f(y_k),\\
		x_{k+1}&=y_k-\eta \widehat{\nabla} f(y_k),\\
		y_k&=x_k+ \tilde{\alpha} d_k.
	\end{align}	
\end{subequations}
These equations are of the form \eqref{eq:control_disc_stochastic} with $\xi_k=[d_k^\T,x_k^\T]^\T\in \R^{2d}$ and
\begin{equation*}
	A = \left[\begin{matrix}(1- \tilde{\alpha}) \tilde{\beta} I_d & 0\\  \tilde{\alpha} I_d& I_d\end{matrix}\right],\quad
	B = \left[\begin{matrix}-\eta(\gamma-1) I_d\\ -\eta I_d\end{matrix}\right],\quad
	C = \left[\begin{matrix}  \tilde{\alpha} I_d & I_d\end{matrix}\right],\quad
	E = \left[\begin{matrix} 0 & I_d\end{matrix}\right].
\end{equation*}

As in deterministic case, the Kronecker product structure of the matrices $A$,$B$, $C$, $E$ lead us to look for a matrix $P$ of the form $P= \widehat{P}\otimes I_d$ for some $2\times2$ matrix $\widehat{P}$ as in \eqref{eq:P} and to set $a_0=1$. Observe that for $P$ of this form and with the matrices $B,E$ given here we have $(EB)^\T(EB) = \eta^2 I_d$ and $$B^\T PB = \eta^2(p_{11}(\gamma-1)^2+2p_{12}(\gamma-1)+p_{22})I_d.$$
Therefore the conditions \eqref{eq:SGC} hold for any $f$ which satisfies \eqref{eq:SGCf} provided that
\begin{equation}\label{eq:BPB_cond}
	p_{11}(\gamma-1)^2+2p_{12}(\gamma-1)+p_{22} \geq 0.
\end{equation}
Now, by Theorem \ref{theo:main_disc_stochastic}, it remains to find $\widehat{P}$ such that $\widehat{T}\preceq 0$, $\tilde{P}=\widehat{P}+(m/2)\widehat{E}^\T \widehat{E} \succ 0$ and \eqref{eq:BPB_cond} holds for $T$ given by \eqref{eq:LMI_disc_stochastic}. The elements of $\widehat{T}$ are
\begin{align*}
	t_{11}&=  \tilde{\alpha} \left( \tilde{\alpha} p_{22}- \tilde{\beta} p_{12}\left( \tilde{\alpha} -1\right)\right)-p_{11}\rho ^2- \tilde{\beta} \left( \tilde{\alpha}  p_{12}- \tilde{\beta}  p_{11} \left( \tilde{\alpha} -1\right)\right) \left( \tilde{\alpha} -1\right)-\frac{ \tilde{\alpha} ^2 m
	}{2},\\
	t_{12}&= \tilde{\alpha}  p_{22}-p_{12} \rho ^2+\frac{ \tilde{\alpha}  m \left(\rho ^2-1\right)}{2}- \tilde{\beta}  p_{12} \left( \tilde{\alpha} -1\right),\\
	t_{13}&=\frac{ \tilde{\alpha}  \rho ^2}{2}-\eta  \left( \tilde{\alpha}  p_{22}- \tilde{\beta}  p_{12} \left( \tilde{\alpha} -1\right)\right)-\frac{ \tilde{\alpha}  \left(\rho ^2-1\right)}{2}\\
&\qquad\qquad\qquad\qquad-\eta  \left( \tilde{\alpha}  p_{12}- \tilde{\beta}  p_{11} \left( \tilde{\alpha} -1\right)\right) \left(\gamma -1\right),\\
	t_{22}&=-\frac{\left(1-\rho ^2\right) \left(m-2 p_{22}\right)}{2},\\
	t_{23} &=-\frac{\rho ^2}{2}-\eta  p_{22}-\eta  p_{12} \left(\gamma -1\right)+\frac{1}{2},\\
	t_{33}&= \frac{L \eta ^2 \rho _{0}}{2}-\eta+\eta^2  p_{22} \rho _{0}+2\eta^2 p_{12} \rho _{0} \left(\gamma -1\right)+\eta^2  p_{11} \rho _{0} \left(\gamma -1\right)^2.
\end{align*}
Following the same  reasoning as in the deterministic case, we first impose that $t_{13}=t_{23}=0$ by setting
\begin{align*}
	p_{11}&=-\frac{ \tilde{\alpha} +2  \tilde{\alpha}  \eta  p_{12}-2  \tilde{\alpha}  \eta  p_{22}-2  \tilde{\beta}  \eta  p_{12}+2  \tilde{\alpha}   \tilde{\beta}  \eta  p_{12}-2  \tilde{\alpha}  \eta  \gamma  p_{12}}{2  \tilde{\beta}  \eta -2  \tilde{\alpha}   \tilde{\beta}  \eta -2  \tilde{\beta}  \eta  \gamma +2  \tilde{\alpha}   \tilde{\beta}  \eta  \gamma },\\
	p_{12}&=\frac{(1-\rho^2)-2\eta p_{22}}{2 \eta(\gamma-1)}.
\end{align*}

In \cite{VBS19} the parameters are set as follows for $f\in \mathcal{F}_{m,L}$ satisfying \eqref{eq:SGCf}
\begin{align}\label{eq:Vaswani_params}
	 \tilde{\alpha}= \frac{\sqrt{m}}{\sqrt{m}+\rho_0\sqrt{L}},
	\quad \tilde{\beta}  = 1-\frac{\sqrt{m}}{\rho_0\sqrt{L}},
	\quad\gamma  = \frac{\sqrt{L}}{\sqrt{m}},
	\quad\eta = \frac{1}{\rho_0L}.
\end{align}
For this choice of parameters with $p_{11}$ and $p_{22}$ set as above we have that
\begin{equation*}
	t_{33}=-\frac{\left(\rho _{0}-1\right) \left(\frac{\sqrt{m}}{\sqrt{L}\rho_0}-1+\rho ^2\right)}{2\sqrt{L} \rho _{0} L \rho _{0} \left(\sqrt{L} \rho _{0}-\sqrt{m}\right)}.
\end{equation*}
If $\rho_0>1$ then $t_{33}$ is only negative for $\rho^2\leq 1-1/(\sqrt{\kappa}\rho_0)$ which gives the same rate as that obtained in \cite[Theorem~2]{VBS19}. Indeed for this choice of $\rho^2$ and parameters as in \eqref{eq:Vaswani_params} one can show that setting
$$
\widehat{P}=\frac{m}{2}\left[\begin{matrix}
	 1 & 1\\
	  1 &  1
\end{matrix}\right]
$$
leads to $T\preceq 0$ with $\rho^2\leq 1-1/(\sqrt{\kappa}\rho_0)$, \modkz{where $\kappa$ is the condition number of the function $f$}. However by choosing parameter values different from \eqref{eq:Vaswani_params}, it is possible to derive improved rates of convergence. We proceed as follows.
We  solve $t_{33}=0$ in terms of $\gamma$ to find
\begin{equation*}
	\gamma = 1-\frac{2  \tilde{\beta} (1- \tilde{\alpha}) -(1-\rho^2) \tilde{\beta}  \rho _{0}+ \tilde{\alpha}   \tilde{\beta} (1-\rho^2) \rho _{0}-L (1- \tilde{\alpha}) \tilde{\beta}  \eta  \rho _{0}}{ \tilde{\alpha}  \rho ^2 \rho _{0}}.
\end{equation*}
We keep the values $\tilde{\alpha}, \tilde{\beta}$ and $\eta$  in \eqref{eq:Vaswani_params}, which results in
\begin{equation}\label{eq:gammaorig}
	\gamma= \frac{\sqrt{\kappa} \rho _{0}-1-(r-1)\rho _{0}}{\rho _{0}-\kappa^{-1/2}r}.
\end{equation}
\begin{figure}[t]
	\centering
	\begin{subfigure}[b]{0.49\textwidth}
		\centering
		\includegraphics[width=\textwidth]{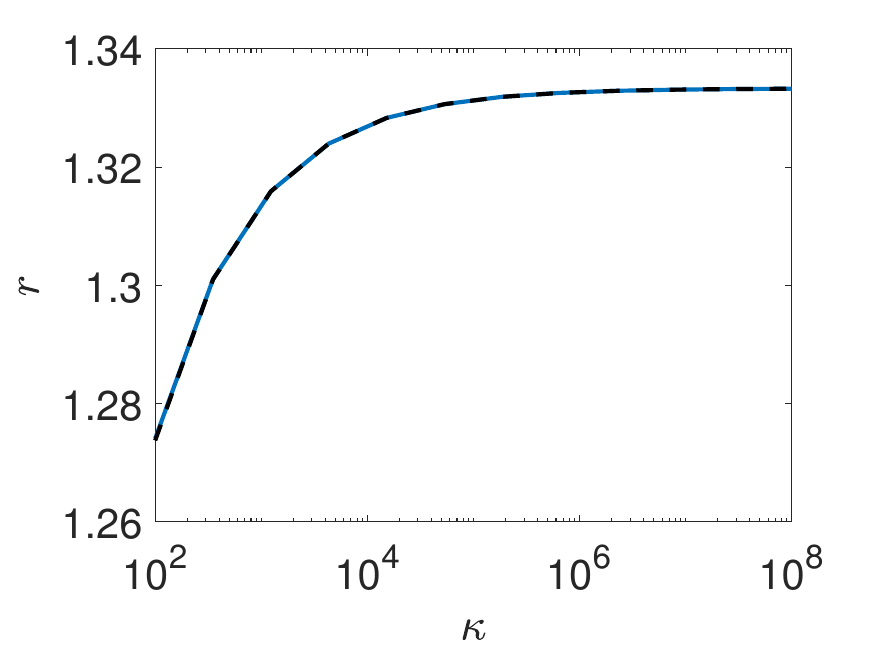}
		\label{fig:r}
	\end{subfigure}
	\hfill
	\begin{subfigure}[b]{0.49\textwidth}
		\centering
		\includegraphics[width=\textwidth]{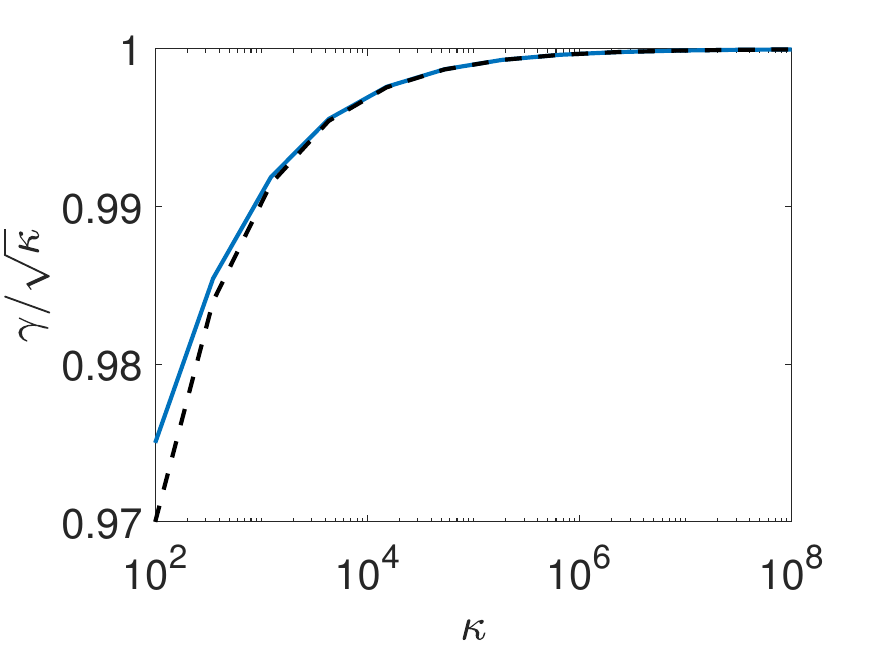}
		\label{fig:gamma}
	\end{subfigure}
%	\begin{subfigure}[b]{0.32\textwidth}
%		\centering
%		\includegraphics[width=\textwidth]{stochastic_constant_curve.eps}
%		\label{fig:multiplitive_constant}
%	\end{subfigure}
	\caption{Convergence of the stochastic Nesterov algorithm with $ \tilde{\alpha},  \tilde{\beta}, \eta$ given by \eqref{eq:Vaswani_params}  when $\rho_0=10$. On the left we have the convergence rate $r$, on the right we show, as a fraction of $\sqrt{\kappa}$, the value of $\gamma$ from \eqref{eq:gammaorig}.  In the dashed lines we show the values when using the approximation $\gamma = \sqrt{\kappa}-(1/3)(1-\rho_0^{-1})$. }
	\label{fig:stochastic_params}
\end{figure}
It remains to consider the $2\times 2$ matrix
\begin{equation*}
	T^0=\left[\begin{matrix}
		t_{11} & t_{12}\\
		t_{12} & t_{22}
	\end{matrix}\right].
\end{equation*}
%Observe that $t_{22}\leq 0$ whenever $p_{22}\leq m/2$.
Following the approach for the deterministic Nesterov algorithm, we calculate $r$ and $p_{22}$ by imposing $\Delta = 0$ and $(\partial/\partial p_{22}) \Delta = 0$ where $\Delta=\det(T^0)$. Since $\Delta$ is a quadratic function of $p_{22}$, there is a unique value of $p_{22}$ which solves $(\partial/\partial p_{22}) \Delta = 0$. Then we have $p_{22}$ as a function of $\rho^2$ and it remains to solve $\Delta=\det(T^0)$ for $\rho^2$ and check the following conditions:
\begin{subequations}
	\begin{align}
		& t_{11}\leq 0 \text{ and } t_{22}\leq 0,\\
		& \widehat{P}+\frac{m}{2}\widehat{E}^\T \widehat{E} \succ 0,\\
		& p_{11}(\gamma-1)^2+2p_{12}(\gamma-1)+p_{22} \geq 0. \label{eq:SGC_check}
	\end{align}
\end{subequations}
The first of these along with having $t_{13}=t_{23}=t_{33}$ and $\Delta=0$ ensures that $T\preceq 0$. The second condition is an assumption in Theorem \ref{theo:main_disc_stochastic} which ensures that the Lyapunov function used upper bounds the Euclidean norm. The third condition is used to ensure that \eqref{eq:SGC} holds.

\begin{figure}[t]
	\centering
	\begin{subfigure}[b]{0.49\textwidth}
		\centering
		\includegraphics[width=\textwidth]{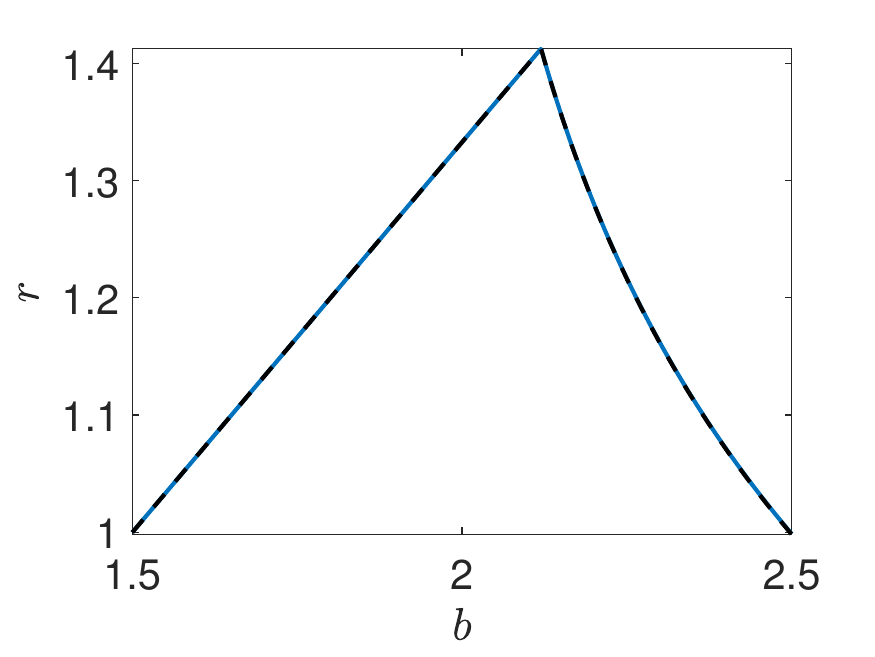}
		\label{fig:r_b}
	\end{subfigure}
	\hfill
	\begin{subfigure}[b]{0.49\textwidth}
		\centering
		\includegraphics[width=\textwidth]{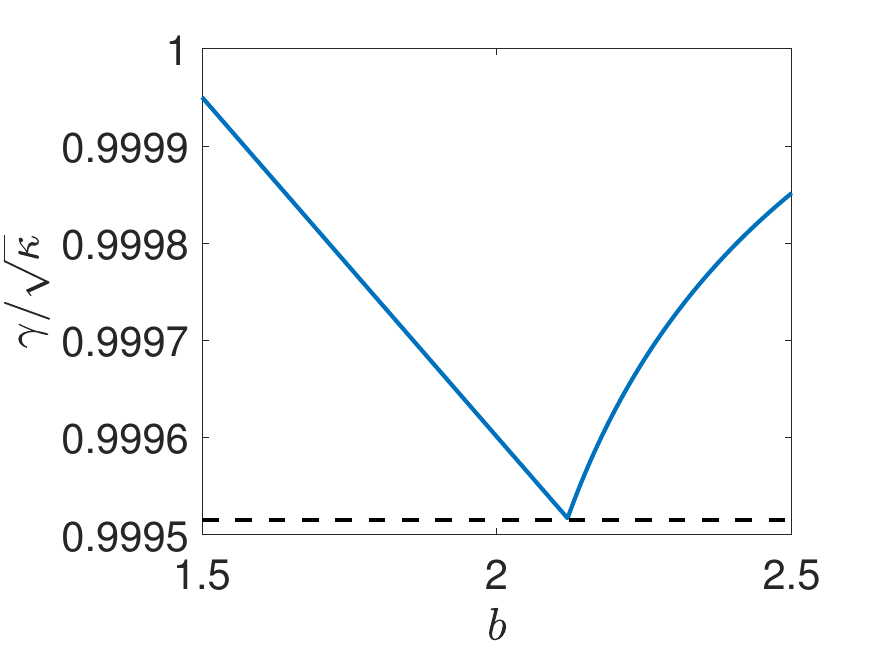}
		\label{fig:gamma_b}
	\end{subfigure}
	%	\begin{subfigure}[b]{0.32\textwidth}
		%		\centering
		%		\includegraphics[width=\textwidth]{stochastic_constant_curve.eps}
		%		\label{fig:multiplitive_constant}
		%	\end{subfigure}
	\caption{Convergence of the stochastic Nesterov algorithm with $ \tilde{\alpha},  \tilde{\beta}, \eta$ given by \eqref{eq:params_with_b}, $\kappa=10^6$ and $\rho_0=10$. On the left we have the convergence rate $r$, on the right we show the optimal choice of $\gamma$ as a fraction of $\sqrt{\kappa}$. In the solid line we show the value of $r$, and the value of $\gamma$ resulting from $t_{33}= 0$. In the dashed lines we show the values when $\gamma$ is set by \eqref{eq:gamma_approx_b}. }
	\label{fig:stochastic_params_b}
\end{figure}

It is convenient to express the variable $\rho$ determined by the procedure above in terms of a new variable $r$ as follows
\begin{equation*}
	\rho^2 = 1-r\frac{\sqrt{m}}{\rho_0\sqrt{L}}.
\end{equation*}
Note that $r=1$ corresponds to the rate obtained in \cite{VBS19} and that therefore values $r>1$ indicate an improved rate.
In Figure \ref{fig:stochastic_params} we show how  $r$ varies as a function of $\kappa$ along with the associated value of $\gamma$ from \eqref{eq:gammaorig}. We see, for $\kappa$ large, $r$ converges to $4/3$ and hence $\gamma$ is approximately
\begin{equation*}
	\gamma \approx \frac{\sqrt{\kappa} -\rho _{0}^{-1}-\frac{1}{3}}{1-\frac{4}{3}\rho _{0}^{-1}\kappa^{-\frac{1}{2}}} \approx \sqrt{\kappa}-\frac{1}{3}(1-\rho_0^{-1}).
\end{equation*}
In the dashed line of Figure \ref{fig:stochastic_params} we show the value of $r$ obtained if we use the approximation of $\gamma= \sqrt{\kappa}-({1}/{3})(1-\rho_0^{-1})$. We see that for all values of $\kappa$ considered we have $r>1$ and for large values of $\kappa$ that $r$ approaches $4/3$.

To leading order in $\kappa$ we have that $P$ matches the matrix $\widehat{\bar{P}}$ in the continuous deterministic setting, indeed
\begin{equation*}
	\widehat{P}\approx \frac{m}{2}
		\left[\begin{array}{cc} 1 & r\\ r & {r^2}/{2} \end{array}\right];
\end{equation*}
from this we see for $\kappa$ sufficiently large that \eqref{eq:SGC_check} holds.

%Here we have established the geometric rate of decay, let us finish by discussing the constant term. We will convert the convergence estimate \eqref{eq:conv_disc_stoch}  into an estimate from the Euclidean norm to itself by using that $f$ is $L$-smooth to establish
%\begin{equation*}
%	\lVert x_k-x^{\star}\rVert^2 \leq K \lVert x_0-x^{\star}\rVert^2 \rho^{2k}
%\end{equation*}

%\begin{equation}\label{eq:multiplicative_constant}
%	K=\frac{\max \sigma (P+\frac{L}{2}E^\T E)}{\min\sigma (\tilde{P})}.
%\end{equation}
%For $L$ large we have $K\approx 25.6\kappa$.

\begin{remark}
	In the preceding analyis, we have chosen to use $ \tilde{\alpha}, \tilde{\beta}$ and $\eta$ as in \cite{VBS19} and set an alternative value for $\gamma$. When $\rho_0=1$, \eqref{eq:Nesterov_stoch_sys} is the same algorithm as \eqref{eq:nest3} except with a different set of parameters, and using the parameters given by \eqref{eq:Vaswani_params} corresponds to setting $\beta=(\sqrt{\kappa}-1)/(\sqrt{\kappa}+1)$. As discussed in Section \ref{sec:polyak}, the choice $\beta=1-\sqrt{2/\kappa}$ allows to show an improved convergence rate. We obtain analogous behaviour to Section \ref{sec:polyak} by proceeding as above but using the parameter choice:
	\begin{equation}\label{eq:params_with_b}
		\tilde{\alpha} = \frac{1}{\rho_0\sqrt{\kappa}},\quad \tilde{\beta} = \frac{1-b\frac{1}{\rho_0\sqrt{\kappa}}}{1-\tilde{\alpha}},\quad \eta = \frac{1}{L\rho_0}.
	\end{equation}
	\modk{Here $b$ is a new parameter to be determined which is introduced to be analogous to the parameter $b$ in Section \ref{sec:polyak}, the choice $b=2$ corresponds to the parameter choice \eqref{eq:Vaswani_params} to leading order.}
As before, we can either use $\gamma$ obtained by solving $t_{33}=0$ for $\gamma$ or an approximation, which now is given by
	\begin{equation}\label{eq:gamma_approx_b}
		\gamma =  \frac{3\,\kappa^{-1/2}+\sqrt{\kappa}\,\rho_0 -\sqrt{2}\,\kappa^{-1/2}-\frac{3\,\sqrt{2}}{2}+\rho_0 -\sqrt{2}\rho_0 }{\rho_0 -\sqrt{2}\kappa^{-1/2}}.
	\end{equation}
	By the same strategy as above, we establish the convergence of \modkz{$\mathbb{E}\norm{x_{k}-x^\star}^2$} with rate $\rho^2$, for $\rho^2=1-r\kappa^{-1/2}\rho_0^{-1}$. In Figure  \ref{fig:stochastic_params_b} we show how $r$ depends on $b$. As in Figure \ref{fig:conta}, we see that  $b=3\sqrt{2}/2$  gives $r=\sqrt{2}$ to leading order in $\kappa$.
\end{remark}

{\bf Acknowledgements.} PD and KCZ acknowledges support from the EPSRC grant EP/V006177/1. JMS has been funded by Ministerio de Ciencia e Innovaci\'{o}n (Spain), project PID2022-136585NB-C21,  MCIN/AEI/10.13039/501100011033/FEDER, UE.

\bibliographystyle{abbrv}
\bibliography{referencesB}

\appendix

\modp{
\section{Expressing properties of $\mathcal{F}_{m,L}$ as matrix inequalities}
}

\modp{
In \cite{MR97} integral quadratic constraints (IQC) were proposed as a tool to describe classes of nonlinearities in control theory. This was adapted for optimization algorithms in \cite{LRP16}. As explained in Section \ref{subsec:MatrixIneq} the key idea here is to express concepts like $m$-strong convexity in terms of matrix inequalities. From such inequalities we  obtain the following lemma in \cite{FRMP18}.}

\modp{\begin{lemma}{\cite[Lemma 4.1]{FRMP18}}
	Fix $f\in\mathcal{F}_{m,L}$. Define $e_k=((\xi_k-\xi^\star)^\T,u_k^\T)^\T$ then the following inequalities hold for all $k$
	\begin{align}
		&e_k^\T M^{(1)}e_k\geq f(x_{k+1})-f(x_k) ,\label{eq:M1ineq}\\
		&e_k^\T M^{(2)}e_k \geq f(x_{k+1})-f(x^\star), \label{eq:M2ineq}\\
		&e_k^\T M^{(3)}e_k\geq 0, \label{eq:M3ineq}
	\end{align}
	where $M^{(1)}$, $M^{(2)}$, $M^{(3)}$ are as in Theorem~\ref{theo:main_disc}.	
\end{lemma}
}

\subsection{Proof of Theorem~\ref{theo:main_con}}\label{subsec:proof_main_con}

\begin{proof}
	\modp{Fix $f\in \mathcal{F}_{m,L}$, and define $V$ by \eqref{eq:cont_lyap}. The estimate \eqref{eq:conv_cont} follows if we establish the following two properties:
		\begin{enumerate}
			\item The function $V(\xi(t),t)$ is non-increasing in $t$.
			\item We can bound $V(\xi,t)$ from below with $\lVert x(t)-x^\star\rVert^2$, i.e.
			\begin{equation}\label{eq:lower_bound_V_cont}
				\lVert x(t)-x^\star\rVert^2 \leq e^{-\lambda t}\frac{\max\sigma(\bar{C}^\T \bar{C})}{\min \sigma(\tilde{P})} V(\xi(t),t).
			\end{equation}
		\end{enumerate}
	}
	
	\modp{
		By means of \eqref{eq:upper_bound}, we can bound $V(\xi(t),t)$ from below with
		\begin{equation*}
			V(\xi(t),t)\geq e^{\lambda t} \left(\frac{m}{2}\lVert x(t)-x^\star\rVert^2+(\xi(t)-\xi^\star)^{\T}
			\bar{P}(\xi(t)-\xi^\star) \right).
		\end{equation*}
		Writing $\tilde{P}=\frac{1}{2}m\bar{C}^\T \bar{C}+\bar{P}$ and using that $x(t)-x^\star=C(\xi(t)-\xi^\star)$, we have
		\begin{equation*}
			V(\xi(t),t)\geq e^{\lambda t} (\xi(t)-\xi^\star)^{\T}
			\tilde{P}(\xi(t)-\xi^\star).
		\end{equation*}
		By the assumptions of Theorem~\ref{theo:main_con} we have that $\tilde{P}\succ 0$ so the quadratic form $(\xi(t)-\xi^\star)^{\T}
		\tilde{P}(\xi(t)-\xi^\star)$ is non-degenerate and bounded from below in terms of the minimum eigenvalue of $\tilde{P}$:
		\begin{equation*}
			V(\xi(t),t)\geq \min\sigma(\tilde{P})e^{\lambda t} \lVert\xi(t)-\xi^\star\rVert^2.
		\end{equation*}
		Using that $x(t)-x^\star=\bar{C}(\xi(t)-\xi^\star)$ we obtain the bound \eqref{eq:lower_bound_V_cont}.
	}
	
	\modp{It remains to show that $V(\xi(t),t)$ is non-increasing in $t$, i.e.\ $\dot{V}(\xi(t),t)\leq 0$. Using that $x=\bar{C}\xi$, we can differentiate $V$ to find
		\begin{eqnarray*}
			\dot{V}(\xi,t)&=&e^{\lambda t} \left(\lambda(f(x)-f(x^{\star}))+\lambda (\xi-\xi^\star)^{\T}
			\bar{P}(\xi-\xi^\star))\right.\\&&\left.\qquad\qquad+\nabla f(x)^\T \bar{C}\dot{\xi}+\dot{\xi}^{\T}
			\bar{P}(\xi-\xi^\star) +(\xi-\xi^\star)^{\T}
			\bar{P}\dot{\xi}\right).
		\end{eqnarray*}
		Replacing $\nabla f(x)$ with $u$ and $\dot{\xi}$ from  \eqref{eq:con_system}
		\begin{align*}
			e^{-\lambda t}\dot{V}(\xi,t)&= \lambda(f(x)-f(x^{\star}))+\lambda (\xi-\xi^\star)^{\T}
			\bar{P}(\xi-\xi^\star))+u^\T \bar{C}(\bar{A}\xi+\bar{B}u)\\
			&\qquad\qquad+(\bar{A}\xi+\bar{B}u)^{\T}
			\bar{P}(\xi-\xi^\star) +(\xi-\xi^\star)^{\T}
			\bar{P}(\bar{A}\xi+\bar{B}u).
		\end{align*}
		Define  $e(t)=[(\xi(t)-\xi^\star)^\T,u(t)^\T]^\T$. Then we can write
		\begin{align*}
			e^{-\lambda t}\dot{V}(\xi,t)&=  \lambda(f(x)-f(x^{\star})) +e(t)^\T\left[\begin{matrix}
				\lambda \bar{P} + \bar{A}^\T \bar{P}+\bar{P}\bar{A}&  \frac{1}{2}(\bar{C}\bar{A})^\T+\bar{P}\bar{B}\\
				\frac{1}{2}\bar{C}\bar{A} +\bar{B}^\T\bar{P} & \frac{1}{2}(\bar{C}\bar{B}+(\bar{C}\bar{B})^\T)
			\end{matrix}\right]e(t),
		\end{align*}
		and
		\begin{align*}
			e^{-\lambda t}\dot{V}(\xi,t)&=  \lambda(f(x)-f(x^{\star})) +e(t)^\T(M^{(0)}+M^{(1)})e(t).
		\end{align*}
	}
	
	\modp{
		In order to control the $f$ terms, we use the following convexity inequality
		\begin{equation}\label{eq:strong_convex2}
			f(y_1)-f(y_2)\leq \nabla f(y_1)^\T (y_1-y_2)-\frac{m}{2}\lVert y_1-y_2\rVert^2.
		\end{equation}
		%Note that setting $y_1=x_{k+1}$ and $y_2=x^\star$ in \eqref{eq:strong_convex2} corresponds to \eqref{eq:M2ineq} in the discrete setting, similarly we can rewrite \eqref{eq:strong_convex2} with $y_1=y$ and $y_2=y^\star$ as
		We rewrite \eqref{eq:strong_convex2} with $y_1=x$ and $y_2=x^\star$ as
		\begin{eqnarray*}
			f(x)-f(x^\star)&\leq& \left[\begin{matrix}
				\xi-\xi^\star\\
				u
			\end{matrix}\right]^\T \left[\begin{matrix}
				\bar{C} & 0_d\\
				0_d & I_d
			\end{matrix}\right]^\T\left[\begin{matrix}
				-\frac{m}{2}I_d & \frac{1}{2}I_d\\
				\frac{1}{2}I_d & 0_d
			\end{matrix}\right]\left[\begin{matrix}
				\bar{C} & 0_d\\
				0_d & I_d
			\end{matrix}\right]\left[\begin{matrix}
				\xi-\xi^\star\\
				u
			\end{matrix}\right]\\ & = &e(t)^\T M^{(2)}e(t).
		\end{eqnarray*}
		Therefore,
		\begin{align*}
			e^{-\lambda t}\dot{V}(\xi,t)& \leq e(t)^\T(M^{(0)}+M^{(1)}+\lambda  M^{(2)})e(t).
		\end{align*}
		By  \eqref{eq:M3ineq}  we have that $M^{(3)} \succeq 0$ and thus, for $\sigma>0$,
		\begin{align*}
			e^{-\lambda t}\dot{V}(\xi,t)& \leq e(t)^\T(M^{(0)}+M^{(1)}+\lambda  M^{(2)}+\sigma M^{(3)})e(t) = e(t)^\T Te(t).
		\end{align*}
		This shows that $\dot{V}(\xi,t)\leq 0$.
	}
\end{proof}

%\section{Proof of the Theorems}
\subsection{Proof of Theorem~\ref{theo:main_disc}}\label{subsec:proof_main_disc}
\begin{proof}
\modp{Fix $f\in\mathcal{F}_{m,L}$ and define $V_k$ by \eqref{eq:liap_disc}.  In order to show that $V$ is a Lyapunov function we need to establish that $V_k$ is non-increasing along trajectories of the algorithm, i.e. $V_{k+1}(\xi_{k+1})\leq V_k(\xi_k)$, and that $V_k(\xi_k)$ controls the convergence of $\lVert x_k-x^\star\rVert$. Consider the second of these issues, which we settle by showing that $V_k$ is a suitable upper bound for the distance between $x_k$ and $x^\star$.}

\modp{
Using \eqref{eq:upper_bound}, we can bound $V_k(\xi)$ from below by
\begin{equation*}
	V_{k}(\xi_k)\geq \rho^{-2 k} \left(a_{0}\frac{m}{2}\lVert x_k-x^\star\rVert^2+(\xi_k-\xi^\star)^{\T}
	P(\xi-\xi^\star) \right).
\end{equation*}
Writing $\tilde{P}=\frac{1}{2}a_0mE^\T E+P$ and using that $x_k-x^\star=E(\xi_k-\xi^\star)$, we have
\begin{equation*}
	V_{k}(\xi_k)\geq \rho^{-2 k} (\xi_k-\xi^\star)^{\T}
	\tilde{P}(\xi_k-\xi^\star).
\end{equation*}
The quadratic form $(\xi-\xi^\star)^{\T}
\tilde{P}(\xi-\xi^\star)$ is non-degenerate and is bounded from below in terms of  the minimum eigenvalue of $\tilde{P}$:
\begin{equation*}
	V_{k}(\xi_k)\geq \min\sigma(\tilde{P})\rho^{-2 k} \lVert\xi_k-\xi^\star\rVert^2.
\end{equation*}
Using that $x_k-x^\star=E(\xi_k-\xi^\star)$ we obtain the following bound:
\begin{equation*}
	\lVert x_k-x^\star\rVert^2\leq \frac{\max\sigma(E^\T E)}{\min\sigma(\tilde{P})}\rho^{2k}V_k(\xi_k).
\end{equation*}}

\modp{
If $V_{k+1}(\xi_{k+1}) \leq V_k(\xi_k)$, then we will have the desired bound
\begin{equation}\label{eq:lowerbound_disc}
	\lVert x_k-x^\star\rVert^2\leq \frac{\max\sigma(E^\T E)}{\min\sigma(\tilde{P})}\rho^{2k}V_0(\xi_0).
\end{equation}
}

\modp{We establish that $V_k(\xi_k)$ is non-increasing as follows. We have
\begin{align}
	&e_k^\T M^{(0)}e_k=(\xi_{k+1}-\xi_{k+1}^\star)^{\T}
	P(\xi_{k+1}-\xi^\star)-\rho^2(\xi_{k+1}-\xi_k^\star)^{\T}
	P(\xi_k-\xi^\star),\label{eq:M0ineq},
\end{align}
where $e_k=((\xi_k-\xi^\star)^\T , u_k^\T )^\T .$
Indeed, by using \eqref{eq:control_disc}, we may write
\begin{align*}
	(\xi_{k+1}-\xi_{k+1}^\star)^{\T}
	P(\xi_{k+1}-\xi^\star) &= (\xi_{k}-\xi_{k}^\star)^{\T}A^\T
	PA(\xi_{k}-\xi^\star)+(u_{k})^{\T}B^\T
	PA(\xi_{k}-\xi^\star)\\
	&\qquad\qquad+(\xi_{k}-\xi_{k}^\star)^{\T}A^\T
	PBu_{k}+u_k^{\T}B^\T
	PBu_k \\
	&=e_k^\T M^{(0)} e_k+\rho^2(\xi_{k+1}-\xi_k^\star)^{\T}
	P(\xi_k-\xi^\star).
\end{align*}
Using \eqref{eq:M1ineq}, \eqref{eq:M2ineq}, and \eqref{eq:M0ineq},
\begin{align*}
	V_{k+1}(\xi_{k+1}) -V_k(\xi_k)&= a_0\rho^{-2(k+1)} \Big(\rho^2((f(E\xi_{k+1})-f(E\xi_k))\\
&\qquad.+(1-\rho^2)(f(E\xi_{k+1})-f(x^\star))\Big)\\
	&\qquad+\rho^{-2(k+1)}\Big((\xi_{k+1}-\xi_{k+1}^\star)^{\T}
	P(\xi_{k+1}-\xi^\star)\\
&\qquad-\rho^2(\xi_{k}-\xi^\star)^{\T}
	P(\xi_{k}-\xi^\star)\Big)\\
	&\leq \rho^{-2(k+1)} e_k^\T \left(M^{(0)}+a_0\rho^2M^{(1)}+a_0(1-\rho^2)M^{(2)}\right)e_k.
\end{align*}
}
\end{proof}

Since $e_k^\T M^{(3)}e_k\geq 0$ by \eqref{eq:M3ineq}, we can add $\ell e_k^\T M^{(3)}e_k$, to bound $V_{k+1}(\xi_{k+1})$ $ -V_k(\xi_k)$  by
\begin{align*}
	&\leq \rho^{-2(k+1)} e_k^\T \left(M^{(0)}+a_0\rho^2M^{(1)}+a_0(1-\rho^2)M^{(2)}+\ell M^{(3)}\right)e_k\\ &= \rho^{-2(k+1)}e_k^\T Te_k.
\end{align*}
Therefore, $V_{k+1}(\xi_{k+1})\leq V_k(\xi_k)$.

\subsection{Proof of Theorem~\ref{theo:main_disc_stochastic}}\label{subsec:proof_main_disc_stochastic}

\modp{\begin{proof}
	The proof of this theorem follows the same argument as the proof of Theorem \ref{theo:main_disc}, except for the derivation of $N^{(1)}$ and $M^{(0)}$; therefore we only show how these terms differ.
	%	Since $f$ is $L$-smooth we have the following matrix inequality
	%	\begin{equation}\label{eq:Lipschitz_MI}
		%		f(x_{k+1})-f(y_k) \leq \left[\begin{matrix}x_{k+1}-y_k\\\nabla f(y_k)\end{matrix}\right]^\T \left[\begin{matrix}\frac{L}{2}I_d & \frac{1}{2}I_d\\ \frac{1}{2}I_d & 0\end{matrix}\right]\left[\begin{matrix}x_{k+1}-y_k\\\nabla f(y_k)\end{matrix}\right].
		%	\end{equation}
	Using  \eqref{eq:control_disc_stochastic},
	\begin{equation}\label{eq:express_nabla_f}
		\left[\begin{matrix}x_{k+1}-y_k\\\nabla f(y_k)\end{matrix}\right] = \left[\begin{matrix}EA-C & EB & 0\\0 & 0 & I_d\end{matrix}\right]\left[\begin{matrix}\xi_k-\xi^\star\\\tilde{u}_k\\ \nabla f(y_k)\end{matrix}\right]
	\end{equation}
	and substituting \eqref{eq:express_nabla_f} into \eqref{eq:L_smooth} (with $x=x_{k+1}$ and $y=y_k$),  $f(x_{k+1})-f(y_k)$  may be upperbounded by
	\begin{align*}
		%& \leq \left[\begin{matrix}\xi_k-\xi^\star\\\nabla f(y_k,z_k)\\ \nabla f(y_k)\end{matrix}\right]^\T\left[\begin{matrix}(EA-C)^\T & 0 \\(EB)^\T & 0\\0  & I_d\end{matrix}\right] \left[\begin{matrix}\frac{L}{2}I_d & \frac{1}{2}I_d\\ \frac{1}{2}I_d & 0\end{matrix}\right]\left[\begin{matrix}EA-C & EB & 0\\0 & 0 & I_d\end{matrix}\right]\left[\begin{matrix}\xi_k-\xi^\star\\\nabla f(y_k,z_k)\\ \nabla f(y_k)\end{matrix}\right]\\
		&  \left[\begin{matrix}\xi_k-\xi^\star\\\tilde{u}_k\\ \nabla f(y_k)\end{matrix}\right]^\T\left[\begin{matrix}\frac{L}{2}(EA-C)^\T(EA-C) & \frac{L}{2}(EA-C)^\T EB & \frac{1}{2}(EA-C)^\T \\ \frac{L}{2}(EB)^\T(EA-C) & \frac{L}{2}(EB)^\T EB & \frac{1}{2}(EB)^\T\\ \frac{1}{2}(EA-C) &\frac{1}{2}EB & 0\end{matrix}\right]\left[\begin{matrix}\xi_k-\xi^\star\\\tilde{u}_k\\ \nabla f(y_k)\end{matrix}\right]\\
		&=:\left[\begin{matrix}\xi_k-\xi^\star\\\tilde{u}_k \\ \nabla f(y_k)\end{matrix}\right]^\T\tilde{N}^1\left[\begin{matrix}\xi_k-\xi^\star\\\tilde{u}_k\\ \nabla f(y_k)\end{matrix}\right].
	\end{align*}
	We can expand this matrix inequality as
	\begin{align*}
		f(x_{k+1})-f(y_k) &\leq (\xi_k-\xi^\star)^\T\tilde{N}^1_{11}(\xi_k-\xi^\star)+2(\xi_k-\xi^\star)^\T\tilde{N}^1_{12}\tilde{u}_k\\ &\qquad+2(\xi_k-\xi^\star)^\T\tilde{N}^1_{13}\nabla f(y_k)
		+\tilde{u}_k^\T\tilde{N}^1_{22}\tilde{u}_k\\
&\qquad +2(\nabla f(y_k))^\T\tilde{N}^1_{23}\tilde{u}_k +(\nabla f(y_k))^\T\tilde{N}^1_{33}\nabla f(y_k) .
	\end{align*}
	Taking expectation (conditional on $\xi_k$) using that $\mathbb{E}[\tilde{u}_k|\xi_k]=\nabla f(y_k)$ and \eqref{eq:SGC}, we have
	\begin{align*}
		\mathbb{E}[f(x_{k+1})-f(y_k)|\xi_k] &\leq (\xi_k-\xi^\star)^\T\tilde{N}^1_{11}(\xi_k-\xi^\star)+2(\xi_k-\xi^\star)^\T(\tilde{N}^1_{12}+\tilde{N}^1_{13})\nabla f(y_k)  \\
		&\qquad\qquad+(\nabla f(y_k))^\T(\rho_0\tilde{N}^1_{22}+2\tilde{N}^1_{23}+\tilde{N}^1_{33})\nabla f(y_k).
	\end{align*}
	We can re-express this as a matrix inequality as
	\begin{equation*}
		\mathbb{E}[f(x_{k+1})-f(y_k)|\xi_k] \leq \left[\begin{matrix}\xi_k-\xi^\star\\\nabla f(y_k)\end{matrix}\right]^\T N^{(1)}\left[\begin{matrix}\xi_k-\xi^\star\\ \nabla f(y_k)\end{matrix}\right],
	\end{equation*}
	with $N^{(1)}$ as given in the statement of the theorem.
	%	\begin{equation*}
		%		N^1 = \left[\begin{matrix}\frac{L}{2}(EA-C)^\T(EA-C) & \frac{1}{2}(EA-C)^\T(LEB+1) \\ \frac{1}{2}(LEB+1)^\T(EA-C) & \frac{L\rho_0}{2}(EB)^\T EB+\frac{1}{2}(EB+(EB)^\T)\end{matrix}\right]
		%	\end{equation*}
	%	As in the non-stochastic setting we also have
	%	\begin{equation*}
		%		f(y_k)-f(x_k) \leq e_k^\T N^2e_k.
		%	\end{equation*}
	%	Then adding * and * gives
	%	\begin{equation*}
		%		f(x_{k+1})-f(x_k) \leq e_k^\T(N^1+N^2)e_k = e_k^\T M^1e_k.
		%	\end{equation*}
	%	Similarly
	%	\begin{equation*}
		%		f(x_{k+1})-f(x_\star) \leq e_k^\T(N^1+N^3)e_k = e_k^\TM^2e_k.
		%	\end{equation*}
	The other term to consider is
	\begin{eqnarray*}
		&&\rho^{-2(k+1)}(\xi_{k+1}-\xi^\star)^\T P(\xi_{k+1}-\xi^\star)-\rho^{-2k}(\xi_{k}-\xi^\star)^\T P(\xi_{k}-\xi^\star)   \\&&\qquad\qquad\qquad\qquad\qquad\quad= \rho^{-2(k+1)}\left[\begin{matrix}\xi_k-\xi^\star\\\tilde{u}_k\end{matrix}\right]^\T\left[\begin{matrix} A^\T PA-\rho^2P & A^\T PB\\ B^\T PA & B^\T PB\end{matrix}\right]\left[\begin{matrix}\xi_k-\xi^\star\\\tilde{u}_k\end{matrix}\right]
	\end{eqnarray*}
	Taking expectation, using that $\mathbb{E}[\tilde{u}_k|\xi_k]=\nabla f(y_k)$ and  \eqref{eq:SGC}, we conclude that
	\begin{align*}
		&\mathbb{E}[\rho^{-2(k+1)}(\xi_{k+1}-\xi^\star)^\T P(\xi_{k+1}-\xi^\star)-\rho^{-2k}(\xi_{k}-\xi^\star)^\T P(\xi_{k}-\xi^\star)|\xi_k]
\\ & \qquad\qquad = \rho^{-2(k+1)}\left[\begin{matrix}\xi_k-\xi^\star\\\nabla f(y_k)\end{matrix}\right]^\T M^{(0)}\left[\begin{matrix}\xi_k-\xi^\star\\\nabla f(y_k)\end{matrix}\right].
	\end{align*}
	The remainder of the proof follows the same argument as Theorem \ref{theo:main_disc}.
\end{proof}
}

\end{document}